\documentclass[reprint,twocolumn,superscriptaddress,showkeys,
nofootinbib,notitlepage,amsmath,amssymb,floatfix,preprintnumbers]{revtex4-1}

\usepackage{caption,subcaption}
\usepackage{graphicx}
\usepackage{amsmath,amssymb,float}
\usepackage{natbib}
\usepackage{color}
\usepackage{soul}
\soulregister\cite7
\soulregister\ref7
\soulregister\pageref7
\usepackage{bm}
\usepackage{array}
\usepackage{hyperref}
\usepackage[labelformat=simple]{subcaption}                                           %enables subfigure 
\begin{document}
\title{Scale Invariance of the Homentropic Inviscid Euler Equations with Application to the Noh Problem}
\author{Jesse F.~Giron}
\email{jgiron@lanl.gov}
\affiliation{Applied Physics, Los Alamos National Laboratory, P.O. Box 1663, MS T082, Los Alamos, New Mexico 87545}
\affiliation{Department of Physics, Box 871504, Arizona State University, Tempe, Arizona  85287-1504}
\author{Scott  D.~Ramsey}
\email{ramsey@lanl.gov}
\affiliation{Applied Physics, Los Alamos National Laboratory, P.O. Box 1663, MS T082, Los Alamos, New Mexico 87545}
\author{Roy S.~Baty}
\email{rbaty@lanl.gov}
\affiliation{Applied Physics, Los Alamos National Laboratory, P.O. Box 1663, MS T082, Los Alamos, New Mexico 87545}
\preprint{LA-UR 19-31605}
\date{Nov. 18, 2019}
\begin{abstract}
We investigate the inviscid compressible flow (Euler) equations constrained by an "isentropic" equation of state (EOS), whose functional form in pressure is an arbitrary function of density alone. Under the aforementioned condition, we interrogate using symmetry methods the scale-invariance of the homentropic inviscid Euler equations. We find that under general conditions, we can reduce the inviscid Euler equations into a system of two coupled ordinary differential equations. To exemplify the utility of these results, we formulate two example scale-invariant, self-similar solutions. The first example includes a shock-free expanding bubble scenario, featuring a modified Tait EOS. The second example features the classical Noh problem, coupled to an arbitrary isentropic EOS. In this case, in order to satisfy the conditions set forth in the classical Noh problem, we find that the solution for the flow is given by a transcendental algebraic equation in the shocked density.
\end{abstract}
\maketitle

\section{Introduction}\label{intro}

The inviscid compressible flow (Euler) equations are a powerful tool in the analytical modeling of shock wave propagation. These equations host a number of well-known, canonical shock solutions, including the Sedov-Taylor-von Neumann blast wave~\cite{sedov,korobeinikov}, the Guderley converging shock~\cite{guderley}, Noh's stagnation shock problem~\cite{noh1987errors} (and more generally, a wide variety of Riemann solutions~\cite{menikoff1989riemann}), and other, less well-known solutions owing to Coggeshall~\cite{coggeshall1991analytic}, among others. These generally semi-analytical and sometimes closed-form solutions are valuable for a variety of purposes, including as code verification test problems~\cite{ramsey2012guderley,burnett2018verification}, intermediate asymptotic entities~\cite{barenblatt1,barenblatt2}, and tools to help understand related but more general or application-specific flows~\cite{motz1979physics,atzeni2004physics}.

Many of the aforementioned solutions of the inviscid Euler equations share two common features: (1) they rely on the assumption of an ideal gas equation of state (EOS) closure law in addition to conservation principles for mass, momentum, and total energy, and (2) they are a direct consequence of invariance of the inviscid Euler equations under various sub-groups of the canonical three-parameter scaling group. Explicit discussion of the latter property is approached to varying degrees by Sedov~\cite{sedov}, Meyer-ter-Vehn and Schalk~\cite{meyer1982selfsimilar}, and Coggeshall~\cite{coggeshall1986lie,coggeshall1991analytic,coggeshall1992group}. 

Indeed, the only means for systematically approaching scale or other invariance of the inviscid Euler equations (as opposed to using dimensional analysis or ad hoc methods) is symmetry analysis, also variously referred to as Lie group or group-theoretic techniques. Following the pioneering considerations of Birkhoff~\cite{birkhoff}, formal symmetry analysis of the invscid Euler equations was carried out by Ovsiannkov~\cite{ovsyannikov1} (see also Holm~\cite{holm1976symmetry}, Hutchens~\cite{hutchens}, and Coggeshall~\cite{coggeshall1991analytic}), from which the existence of numerous scaling and other symmetries was rigorously verified. A powerful outcome of the analysis is the generalization beyond the ideal gas law of the included EOS so as to enable the continued presence of scaling or other transformations, and hence their various dependent canonical solutions. For example, for scaling transformations, the most general admissible EOS ``is ... of the Mie-Gruneisen type"~\cite{hutchens}.

Contained within this class of admissible EOS models are the so-called ``isentropic" (also known as ``adiabatic") EOS models as discussed by (for example) Anderson~\cite{anderson1990modern} or Leveque~\cite{leveque}. An EOS of this form is interpreted as being valid along a flow's isentropes. As a result, the relations of pressure $P = f\left(\rho,S\right) \to P = f\left(\rho\right)$ such that internal energy $I = g\left(\rho\right)$, where $f$ and $g$ are arbitrary functions of the fluid density $\rho$, and $S$ is the fluid entropy. A canonical example of an isentropic EOS is the modified Tait EOS, given by
\begin{equation}\label{tait_eos}
P\left(\rho\right)= B\left[\left(\frac{\rho}{\rho_{\rm{ref}}}\right)^\gamma -1\right],
\end{equation}
for constant $B$, $\gamma$, and $\rho_{\rm{ref}}$, as discussed by Zel'dovich and Raizer~\cite{zeldovich_raizer}. 

The applications of isentropic EOS classes such as the modified Tait form appear in numerous practical contexts. In general, Stanyukovich~\cite{stanyukovich2016unsteady} notes that EOS classes of this form assume ``...a great deal of accuracy" for ``...purely adiabatic processes in substances of low compressibility." Notably included in these scenarios are a variety of processes in water or other liquids~\cite{li1967equation}. For example, both Ridah~\cite{ridah1988shock} and Baum, et al.~\cite{baum1959physics} discuss relevant scenarios including shock propagation in liquids, and liquid propellant problems. Baum, et al.~\cite{baum1959physics}, Wardlaw and Mair~\cite{wardlaw1998spherical}, and Cole~\cite{cole1948underwater} provide extensive treatises on underwater explosion processes, in which the invocation of modified Tait or other "genrealized isentropic" EOS forms play no small part. The assumption of an isentropic EOS is also central to numerous studies associated with bubble collapse, cavitation and cavitation-induced damage, and sonoluminesence ~\cite{akinsete1969nonsimilar,hunter1960collapse}. In general, Stanyukovich~\cite{stanyukovich2016unsteady}, Baum, et al.~\cite{baum1959physics}, and Zel'dovich and Raizer~\cite{zeldovich_raizer} note that isentropic EOS models are useful from the standpoint of intermediate-pressure (e.g., $10$~atm $< P <$ $100,000$~atm or $1.0\times 10^{-5}\rm{ Mbar}$ $< P <$ $1.0\times 10^{-1}\rm{ Mbar}$) scenarios in condensed materials; in addition to processes in liquids, these regimes also include detonating explosives, the impact of detonation products against metallic surfaces, and projectile impacts against the ground.

Moreover, when coupled to the inviscid Euler equations, an isentropic EOS can obviate the need for an energy conservation law; \footnote{As will be shown in Sec.~\ref{Euler_constrained_EOS}, total energy conservation is still included in the inviscid Euler system, but becomes a redundant condition under the assumption of an isentropic EOS.} as such, the structure of the inviscid Euler equations is simpler than their traditional counterparts explicitly featuring a total energy conservation relation. In addition, the assumption of an isentropic EOS also includes homentropic flow (i.e., the fluid entropy is explicitly constant in all space and time) as a sub-case; though the presence of discontinuous shock waves is still possible through the presence of piecewise constant entropy solutions.

Given the broad and important application space of compressible flows featuring isentropic EOS classes, the goal of this work is to provide a dedicated symmetry analysis of the (piecewise) homentropic inviscid Euler equations, coupled to an arbitrary isentropic EOS. The outcomes of this analysis are complementary to those of Ovsiannikov~\cite{ovsyannikov1}, Holm~\cite{holm1976symmetry}, and others, in that conditions for the presence of scale-invariance on the isentropic EOS can be derived and compared to existing results in other contexts. Any such comparison will thus illuminate what symmetries are lost, preserved, or gained in moving between the traditional and homentropic inviscid Euler equation settings.

Furthermore, a second goal of this work is to leverage the symmetry analysis results to construct analogs of classical solutions of the inviscid Euler equations: for example the Noh stagnation shock problem. Some work along these lines has recently been performed by Ramsey et al.~\cite{Ramsey_Boyd_Burnett}, Burnett et al.~\cite{burnett2018verification}, Velikovich et al.~\cite{velikovich2018solution}, and Deschner et al~\cite{deschner2018self}. The motivation for selecting a Noh problem for demonstration purposes is thus to enable comparison to this wide body of existing work.

In support of these goals, Sec.~\ref{mathematical_model} includes a review of the salient mathematical model, including the reduction of the inviscid Euler equations to homentropic form following inclusion of an isentropic EOS, and construction of the associated shock jump conditions. Sec.~\ref{symmetry_analysis} provides a brief review of the differential form or ``isovector" formalism for conducting symmetry analysis of differential equations. Sec.~\ref{scaling_analysis} features a symmetry analysis of the homentropic inviscid Euler equations and shock jump conditions, with an emphasis on scaling phenomena. Sec.~\ref{examples} provides two example scale-invariant solutions of the underlying mathematical model; notably including further symmetry analysis of the conditions unique to the Noh stagnation shock problem, an associated reduction of the homentropic inviscid Euler equations to ODEs, an analysis of the remaining conditions for obtaining a physically relevant solution, and an example solution for the modified Tait EOS. Finally, a summary and recommendation for future work is provided in Sec.~\ref{conclusion}.

%%%%%%%%%%%%%%%%%%%%%%%%%%%%%%%%%%%%%%%%%%%
%				          %
% 	    Mathematical Model		  %
%					  %
%%%%%%%%%%%%%%%%%%%%%%%%%%%%%%%%%%%%%%%%%%%
\section{Mathematical Model}\label{mathematical_model}
%%%%%%%%%%%%%%%%%%%%%%%%%%%%%%%%%%%%%%%%%%%
%				          %
% 	     Euler Equations		  %
%					  %
%%%%%%%%%%%%%%%%%%%%%%%%%%%%%%%%%%%%%%%%%%%
\subsection{Euler Equations}

As shown by numerous authors~\cite{Ramsey_Baty,harlow1971fluid}, we can write the one-dimensional (1D) inviscid compressible flow (Euler) equations as follows:
\begin{equation}\label{cons_mass}
\frac{\partial\rho}{\partial t} + u \frac{\partial\rho}{\partial r} + \rho\left[ \frac{\partial u}{\partial r} + \frac{nu}{r}\right] =0,
\end{equation}
\begin{equation}\label{cons_mom}
\frac{\partial u}{\partial t} + u\frac{\partial u}{\partial r} +\frac{1}{\rho} \frac{\partial P}{\partial r} =0,
\end{equation}
\begin{equation}\label{cons_energy_e_n_p}
\frac{\partial E}{\partial t} + u\frac{\partial E}{\partial r} + \frac{P}{\rho}\left[\frac{\partial u}{\partial r}+\frac{nu}{r}\right]=0,
\end{equation}
where the mass density $\rho\left(r,t\right)$, radial flow velocity $u\left(r,t\right)$, pressure $P\left(r,t\right)$, and total energy per unit mass $E\left(r,t\right)$ are functions of the position coordinate $r$ and time $t$. The space index $n=0$, $1$, or $2$ corresponds to 1D planar, cylindrical, or spherical geometries, respectively. In addition, $E\left(r,t\right)$ is a function of the specific internal energy $I$ (SIE; internal energy per unit mass), that is
\begin{equation}\label{fund_energy}
E\left(r,t\right) = I\left(r,t\right) + \frac{1}{2} \left[u\left(r,t\right)\right]^2. 
\end{equation}
Equations~(\ref{cons_mass})-(\ref{cons_energy_e_n_p}) represent conservation of mass, momentum, and total energy. Using Eqs.~(\ref{cons_mass}), (\ref{cons_mom}), and (\ref{fund_energy}), Eq.~(\ref{cons_energy_e_n_p}) may be rewritten as 
\begin{equation}\label{energy_relation_sub}
\frac{\partial I}{\partial t} + u\frac{\partial I}{\partial r}-\frac{u}{\rho}\frac{\partial P}{\partial r}-\frac{P}{\rho^2}\left[\frac{\partial\rho}{\partial t} + u\frac{\partial\rho}{\partial r}\right] =0,
\end{equation}
and may be further simplified using the fundamental thermodynamic relation~\cite{mandl1988statistical,bowley1999introductory,adkins1983equilibrium,landau2013statistical,zemansky1966basic} between $\rho$, $P$, $I$, the fluid temperature $T$, and the fluid entropy $S$,
\begin{equation}\label{law}
dI= TdS+\frac{P}{\rho^2}d\rho.
\end{equation} 
Using the chain rule in conjunction with Eq.~(\ref{law}), Eq.~(\ref{energy_relation_sub}) becomes
\begin{equation}\label{isentropic_flow}
\frac{\partial S}{\partial t} + u\frac{\partial S}{\partial r} =0, 
\end{equation}
the equation for isentropic flow; this result is expected as dissipative processes (e.g., viscosity and heat conduction) are absent from Eqs.~(\ref{cons_mass})-(\ref{cons_energy_e_n_p}). 
If the entropy $S$ is assumed to be a function of the fluid density $\rho$ and pressure $P$, Eq.~(\ref{isentropic_flow}) may be expanded to yield
\begin{equation}
\frac{\partial S}{\partial \rho}\biggr\rvert_P\left[\frac{\partial \rho}{\partial t}+ u\frac{\partial \rho}{\partial r}\right]+ \frac{\partial S}{\partial P}\biggr\rvert_\rho\left[\frac{\partial P}{\partial t}+ u\frac{\partial P}{\partial r}\right]=0,
\end{equation} 
or substituting Eq.~(\ref{cons_mass})
\begin{equation}\label{cons_energy}
\frac{\partial P}{\partial t} + u\frac{\partial P}{\partial r} + K_S\left[ \frac{\partial u}{\partial r} + \frac{nu}{r}\right] = 0,
\end{equation}
where the adiabatic bulk modulus $K_S\left(\rho,P\right)$ is defined by
\begin{equation}\label{bulk_def}
K_S\left(\rho,P\right)=-\rho\frac{\frac{\partial S}{\partial \rho}\biggr\rvert_P}{\frac{\partial S}{\partial P}\biggr\rvert_\rho},
\end{equation}
or, as shown by Axford~\cite{axford},
\begin{equation}\label{bulkmod}
K_S(\rho,P) = \frac{P}{\rho}\frac{\partial P}{\partial I}\biggr\rvert_\rho + \rho \frac{\partial P}{\partial \rho}\biggr\rvert_I.
\end{equation}
The adiabatic bulk modulus appears only in the total energy (or entropy) conservation relation given by Eq.~(\ref{cons_energy}), and is a measure of the fluid's resistance to uniform, constant entropy compression. It is obtained from an incomplete EOS of the form $P = P\left(\rho,I\right)$ via Eq.~(\ref{bulkmod}), and is also related to the fluid sound speed $c$ by
\begin{equation}\label{sound}
K_S = \rho c^2 .
\end{equation} 
%%%%%%%%%%%%%%%%%%%%%%%%%%%%%%%%%%%%%%%%%%%
%				          %
% 	    Isentropic EOS		  %
%					  %
%%%%%%%%%%%%%%%%%%%%%%%%%%%%%%%%%%%%%%%%%%%
\subsection{Isentropic Equation of State}\label{isentropic_eos}

To further simplify the inviscid Euler equations, we will consider the following ``isentropic" EOS:
\begin{equation}\label{eos}
P = f(\rho),
\end{equation}
where $f$ is an arbitrary function of the fluid density. An example of an isentropic EOS is the modified Tait EOS as discussed in Sec.~\ref{intro}. In this relation, the parameters $B$ and $\gamma$ ostensibly depend on the entropy $S$ but are taken to be material-dependent constants over the pressure regime of interest. 

The adiabatic bulk modulus corresponding to an isentropic EOS can be determined using Eq.~(\ref{bulkmod}):
\begin{equation}\label{abm_isen}
K_S\left(\rho\right) = \rho f'\left(\rho\right),
\end{equation}
where the prime denotes the derivative with respect to the indicated argument. Then, with Eq.~(\ref{bulk_def}), the corresponding fluid entropy may be determined from:
\begin{equation}
\label{entropy_pde_general}
f'\left(\rho\right) \frac{\partial S}{\partial P} + \frac{\partial S}{\partial \rho} = 0,
\end{equation}
which may be solved using the Method of Characteristics to yield:
\begin{equation}
\label{entropy_sol}
S = F\left[P - f\left(\rho\right)\right],
\end{equation}
where $F$ is an arbitrary function of the indicated argument. Given Eq.~(\ref{eos}), any bounded, non-trivial (i.e., giving $S \neq 0$) parameterization of $F$ yields:
\begin{equation}
\label{entropy_def}
S = S_0,
\end{equation}
where $S_0$ is a constant (an example of $F$ that yields Eq.~(\ref{entropy_def}) is $F = S_0$). Since the entropy associated with Eq.~(\ref{eos}) is constant, flows featuring this EOS are referred to as isentropic. For the case where $S$ is also independent of $r$ and $t$ (i.e., constant everywhere, and not just along streamlines; physically, this phenomenon may originate through the judicious selection of initial or boundary conditions), the resulting flows are referred to as homentropic. The homentropic sub-case will be emphasized throughout the remainder of this work.

With this result, the specific internal energy $I$ of the fluid may be calcuated using Eqs.~(\ref{law}) and~(\ref{eos}),
\begin{equation}
\label{sie_def}
dI = \frac{f\left(\rho\right)}{\rho^2} d\rho,
\end{equation}
such that
\begin{equation}
\label{Isolve}
I = g\left(\rho\right) + I_0,
\end{equation}
where $g$ is related to $f$ via
\begin{eqnarray}
g &=& \int \frac{f\left(\rho\right)}{\rho^2} d\rho \label{gdef}, \\
&\rightarrow& f = \rho^2 g' , \label{gpdef}
\end{eqnarray}
and $I_0$ is a constant set by an initial condition (e.g., $I\left(\rho_{\rm{initial}}\right) = 0$, for some reference density $\rho_{\rm{initial}}$).
%
%%%%%%%%%%%%%%%%%%%%%%%%%%%%%%%%%%%%%%%%%%%
%				          %
%        Euler Equation and the EOS       %
%				          %
%%%%%%%%%%%%%%%%%%%%%%%%%%%%%%%%%%%%%%%%%%%
%
\subsection{Reduced Inviscid Euler Equations}\label{Euler_constrained_EOS}
We can first take advantage of Eq.~(\ref{bulkmod}) in order to rewrite Eqs.~(\ref{cons_mom}) and (\ref{cons_energy_e_n_p}), that is,
\begin{eqnarray}
K_S(\rho) &=& \frac{P}{\rho}\frac{\partial P}{\partial I}\biggr\rvert_\rho + \rho \frac{\partial P}{\partial \rho}\biggr\rvert_I,\\
&=&\frac{P}{\rho}\frac{\partial f}{\partial I}\biggr\rvert_\rho + \rho \frac{\partial f}{\partial \rho}\biggr\rvert_I,\\
&=& 0+ \rho \frac{\partial f}{\partial \rho}\biggr\rvert_I,\\
\rightarrow\frac{\partial f}{\partial \rho} \biggr\rvert_I&=& \frac{K_S(\rho)}{\rho}\label{f_of_rho_relation}.
\end{eqnarray}
Substituting Eqs.~(\ref{eos}) and (\ref{f_of_rho_relation}) into Eq.~(\ref{cons_energy}) gives,
\begin{eqnarray}
r\frac{df}{d\rho}\frac{\partial \rho}{\partial t}+ ru\frac{df}{d\rho}\frac{\partial \rho}{\partial r} + r\rho \frac{df}{d\rho}\frac{\partial u}{\partial r} + n\rho u\frac{df}{d\rho}&=&0\label{eos_energy_equation_gone}.
\end{eqnarray}
Dividing by $\frac{df}{d\rho}$ yields Eq.~(\ref{cons_mass}). Therefore, conservation of energy is automatically satisfied and is a redundant condition for an isentropic equation of state. Therefore, the resulting Euler equations collapse to
\begin{eqnarray}
r\frac{\partial\rho}{\partial t} + ru \frac{\partial\rho}{\partial r} + \rho r\frac{\partial u}{\partial r} + n\rho u &=&0,\label{simp_eos_cons_mass}\\
\rho\frac{\partial u}{\partial t} + \rho u\frac{\partial u}{\partial r} +\frac{d f}{d \rho}\frac{\partial\rho}{\partial r}&=&0.\label{eos_cons_mom}
\end{eqnarray}
Substituting Eq. (\ref{f_of_rho_relation}) into Eq. (\ref{eos_cons_mom}) yields: 
\begin{eqnarray}
\rho^2\frac{\partial u}{\partial t} + \rho^2 u\frac{\partial u}{\partial r} +K_S(\rho)\frac{\partial\rho}{\partial r}&=&0.\label{simp_eos_cons_mom}
\end{eqnarray}
As a result, given a specific form of Eq.~(\ref{eos}), Eqs.~(\ref{simp_eos_cons_mass}) and (\ref{simp_eos_cons_mom}) are a closed system of two partial differential equations (PDEs) with respect to $\rho$ and $u$.
%%%%%%%%%%%%%%%%%%%%%%%%%%%%%%%%%%%%%%%%%%%
%				          %
%      Piecewise Homentropic Flows        %
%				          %
%%%%%%%%%%%%%%%%%%%%%%%%%%%%%%%%%%%%%%%%%%%
%
%
\subsection{Piecewise Homentropic Flows}
Like their more general counterparts, the homentropic inviscid Euler equations may admit discontinuous solutions. These solutions are possible provided that mass and momentum are conserved across any discontinuities. This is the case as Eqs.~(\ref{cons_mass})-(\ref{cons_mom}) contain no sources or sinks of these quantities. 

Equations ensuring conservation of mass and momentum across a discontinuity are derived in numerous sources (e.g. Zel'dovich and Raizer~\cite{zeldovich_raizer}). In the current case, they may be written as
\begin{eqnarray}
\left(u_1-D\right)\rho_1 &=& \left(u_2-D\right)\rho_2,\label{rankine_hugoniot_mass_jump}\\
P_1+\rho_1\left(u_1-D\right)u_1 &=& P_2 + \rho_2\left(u_2-D\right)u_2.\label{rankine_hugoniot_mom_jump}
\end{eqnarray}
In Eqs.~(\ref{rankine_hugoniot_mass_jump}) and (\ref{rankine_hugoniot_mom_jump}), the subscripts $1$ and $2$ denote the fluid state immediately adjacent to either side of the discontinuity that propagates with an arbitrary time-dependent velocity $D(t)$. Equations (\ref{rankine_hugoniot_mass_jump}) and (\ref{rankine_hugoniot_mom_jump}) are the general Rankine-Hugoniot discontinuity ``jump conditions" corresponding to Eqs.~(\ref{simp_eos_cons_mass}) and (\ref{simp_eos_cons_mom}). While Eqs.~(\ref{cons_mass}) and (\ref{cons_mom}) are separately valid on either side of the discontinuity, Eqs.~(\ref{simp_eos_cons_mass}) and (\ref{simp_eos_cons_mom}) are essentially internal boundary conditions that are applied at the interface position to connect the regional solutions of Eqs.~(\ref{simp_eos_cons_mass}) and (\ref{simp_eos_cons_mom}) into a global solution.

As was the case with the inviscid Euler equations themselves, the isentropic EOS may be used to further reduce the jump conditions given by Eqs.~(\ref{rankine_hugoniot_mass_jump}) and (\ref{rankine_hugoniot_mom_jump}). In particular, with Eq.~(\ref{eos}), Eqs.~(\ref{rankine_hugoniot_mass_jump}) and (\ref{rankine_hugoniot_mom_jump}) become
\begin{eqnarray}
\left(u_1-D\right)\rho_1 &=& \left(u_2-D\right)\rho_2,\label{rankine_hugoniot_mass_jump_hom}\\
f\left(\rho_1\right)+\rho_1\left(u_1-D\right)u_1 &=& f\left(\rho_2\right) + \rho_2\left(u_2-D\right)u_2.\;\;\;\;\label{rankine_hugoniot_mom_jump_hom}
\end{eqnarray}
These relations represent a general result for an isentropic EOS. Of particular interest to this work are scenarios where Eqs.~(\ref{rankine_hugoniot_mass_jump_hom}) and (\ref{rankine_hugoniot_mom_jump_hom}) characterize the propagation of a shock wave. As discussed by Zel'dovich and Raizer~\cite{zeldovich_raizer}, physical assumptions relevant to shock waves are
\begin{eqnarray}
\rho_2 &>& \rho_1, \\
u_2 &\neq& u_1, \\
P_2 &>& P_1, 
\end{eqnarray}
where the subscripts 1 and 2 denote the unshocked and shocked states, respectively. 

In addition to these relations, Eq.~(\ref{entropy_def}) sets the entropy of a flow associated with an isentropic EOS. As discussed in Sec.~\ref{isentropic_eos}, the assumption of Eq.~(\ref{eos}) results in Eq.~(\ref{entropy_def}). When a shock wave is present, any otherwise homentropic flow is instead piecewise homentropic, with the entropy assuming piecewise constant values on either side of the discontinuity. In this scenario the entropy jump across the shock wave must be strictly positive, as required by the second law of thermodynamics.

%%%%%%%%%%%%%%%%%%%%%%%%%%%%%%%%%%%%%%%%%%%
%				          %
%           Symmetry Analysis             %
%				          %
%%%%%%%%%%%%%%%%%%%%%%%%%%%%%%%%%%%%%%%%%%%
\section{Symmetry Analysis}\label{symmetry_analysis}
%%%%%%%%%%%%%%%%%%%%%%%%%%%%%%%%%%%%%%%%%%%
%				          %
%           Differential Forms            %
%				          %
%%%%%%%%%%%%%%%%%%%%%%%%%%%%%%%%%%%%%%%%%%%
%
%
\subsection{Differential Forms}\label{diff_forms}
A goal of this work is to subject Eqs.~(\ref{simp_eos_cons_mass}) and (\ref{simp_eos_cons_mom}) to symmetry analysis, with an emphasis on invariance under scaling transformations in all variables. In order to affect this procedure, we will employ the ``isovector" approach of Harrison and Estabrook~\cite{harrison_estabrook}, which requires all relevant differential equations be recast as an equivalent exterior differential system (EDS).

For scaling transformations, there are of course several ways to determine invariance of structures such as Eqs.~(\ref{simp_eos_cons_mass}) and (\ref{simp_eos_cons_mom}): direct substitution of global transformations, execution of Lie's~\cite{lie1888theorie1,lie1890theorie2,lie1891vorlesungen,lie1893vorlesungen,lie1893theorie3,lie1896theorie} ``classical" method (see, for example, Ovsiannikov~\cite{ovsyannikov1}, Bluman and collaborators~\cite{bluman_anco,bluman_kumei}, Olver~\cite{olver}, or Cantwell~\cite{cantwell}), or the isovector method (see, for example, Edelen~\cite{edelen2005applied}, Suhubi~\cite{suhubi2013exterior}, and Stephani~\cite{stephani}). This last method possesses several advantages:
\begin{itemize}
\item The isovector method is a more intuitive geometric setting for differential equations.
\item The isovector method obviates the need for sometimes cumbersome prolongation formulae associated with extending infinitesimal group generators to a higher-dimensional manifold. 
\item Using the isovector method, relatively simple symmetry analysis results (e.g., invariance under scaling) are readily extended to analyses featuring more general transformations.
\end{itemize}

The drawback of the isovector formalism is that it requires the relevant differential equations to be equivalently expressible as a first-order system.  This condition is already met in the case of the homentropic inviscid Euler equations. 

To express Eqs.~(\ref{simp_eos_cons_mass}) and (\ref{simp_eos_cons_mom}) as an EDS, they may first be multiplied by the differential volume element $dt\wedge dr$ to yield
\begin{eqnarray}
r d\rho\wedge dr - ru d\rho\wedge dt - \rho r du\wedge dt + n\rho u dt\wedge dr &\equiv& \mu_1,\label{mu_1}\;\;\;\;\;\;\;\;
\\
\rho^2 du\wedge dr -\rho^2 u du \wedge dt - K_S(\rho) d\rho\wedge dt &\equiv& \mu_2\label{mu_2},
\end{eqnarray}
where Eqs.~(\ref{mu_1}) and (\ref{mu_2}) are referred to as a system of ``2-forms"\footnote{A ``0-form" is any function, while a ``1-form" includes only single differentials denoted by $d$. 2-forms are collections of products of two 1-forms, as indicated. Forms of arbitrary order can also in general be defined as necessary.}. In constructing these expressions, it has been implicitly assumed that all partial derivatives may be regarded as quotients of differentials identified by the exterior derivative operator $d$. The operator $\wedge$ used to multiply differentials is known as a wedge product, which has the properties
\begin{eqnarray}
dq_i\wedge dq_j &=& -dq_j\wedge dq_i,\label{antisymmetric}\\
dq_i\wedge dq_i &=& 0,\label{null}
\end{eqnarray}
for all general coordinates $q_i$. More comprehensive overviews of differential geometry are provided in Edelen \cite{edelen2005applied}, Suhubi \cite{suhubi2013exterior}, Bryant et. al. \cite{bryant2013exterior}, and Bourbaki \cite{nicolas1989lie}.  

As written, Eqs.~(\ref{mu_1}) and (\ref{mu_2}) indicate that previously independent and dependent variables are interpreted as entirely independent of each other and represent differential objects in a higher-dimensional manifold. To establish the equivalence between Eqs.~(\ref{mu_1}) and (\ref{mu_2}) and their PDE counterparts Eqs.~(\ref{simp_eos_cons_mass}) and (\ref{simp_eos_cons_mom}), we must enforce the relationship between independent and dependent variables. This process is referred to as ``sectioning'' by Harrison and Estabrook~\cite{harrison_estabrook}. In this process, the solution sub-manifold is chosen by selecting independent and dependent variables; as such, the exterior derivatives of the selected dependent variables become total derivatives in the independent variables:
\begin{eqnarray}
d\rho &=& \frac{\partial \rho}{\partial r}dr + \frac{\partial \rho}{\partial t}dt\label{total_rho_deriv},\\
du &=& \frac{\partial u}{\partial r}dr + \frac{\partial u}{\partial t}dt\label{total_u_deriv}.
\end{eqnarray}
Substituting Eqs.~(\ref{total_rho_deriv}) and (\ref{total_u_deriv}) into Eqs.~(\ref{mu_1}) and (\ref{mu_2}) yields
\begin{eqnarray}
&&r \left[\frac{\partial \rho}{\partial r}dr + \frac{\partial \rho}{\partial t}dt\right]\wedge dr - ru \left[\frac{\partial \rho}{\partial r}dr + \frac{\partial \rho}{\partial t}dt\right]\wedge dt \nonumber\\
&&- \rho r \left[\frac{\partial u}{\partial r}dr + \frac{\partial u}{\partial t}dt\right]\wedge dt + n\rho u dt\wedge dr,
\end{eqnarray}
and
\begin{eqnarray}
&&\rho^2 \left[\frac{\partial u}{\partial r}dr + \frac{\partial u}{\partial t}dt\right]\wedge dr -\rho^2 u \left[\frac{\partial u}{\partial r}dr + \frac{\partial u}{\partial t}dt\right] \wedge dt \nonumber\\
&&- K_S \left[\frac{\partial \rho}{\partial r}dr + \frac{\partial \rho}{\partial t}dt\right]\wedge dt .
\end{eqnarray}
Using the properties from Eqs.~(\ref{antisymmetric}) and (\ref{null}) we have
\begin{eqnarray}
&&\left[r \frac{\partial \rho}{\partial t} + ru \frac{\partial \rho}{\partial r} + \rho r \frac{\partial u}{\partial r}+ n\rho u \right]dt\wedge dr,\\
&&\left[\rho^2 \frac{\partial u}{\partial t}+\rho^2 u \frac{\partial u}{\partial r}+ K_S(\rho) \frac{\partial \rho}{\partial r}\right]dt\wedge dr.
\end{eqnarray}
By setting these relations equal to zero, the nontrivial solution (i.e., $dt\wedge dr \neq 0$) that follows is Eqs.~(\ref{simp_eos_cons_mass}) and (\ref{simp_eos_cons_mom}). This is process is referred to as ``annulling" by Harrison and Estabrook~\cite{harrison_estabrook}.

Given the equivalence between Eqs.~(\ref{simp_eos_cons_mass}) and (\ref{simp_eos_cons_mom}) and Eqs.~(\ref{mu_1}) and (\ref{mu_2}), the latter will be the system of 2-forms subjected to symmetry analysis in the developments to follow.
%%%%%%%%%%%%%%%%%%%%%%%%%%%%%%%%%%%%%%%%%%%
%				          %
%               INVARIANCE                %
%				          %
%%%%%%%%%%%%%%%%%%%%%%%%%%%%%%%%%%%%%%%%%%%
\subsection{Invariance}\label{invariance}

For a continuously variable transformation parameter $\epsilon$ with identity element $\epsilon = 0$, an objective of this work is to determine for what values of the constants $a_1 - a_6$ the global scaling transformations given by 
\begin{eqnarray}
t_{\rm{new}} &=& e^{\epsilon a_1}t,\label{tnew}\\
r_{\rm{new}} &=& e^{\epsilon a_2}r,\label{rnew}\\
\rho_{\rm{new}} &=& e^{\epsilon a_3}\rho,\label{rhonew}\\
u_{\rm{new}} &=& e^{\epsilon a_4}u,\label{unew}\\
P_{\rm{new}} &=& e^{\epsilon a_5}P,\label{Pnew}\\
I_{\rm{new}} &=& e^{\epsilon a_6}I,\label{Inew}
\end{eqnarray}
leave invariant the EDS representation of the homentropic inviscid Euler equations.

The identification of all admissible point-groups associated with the more general inviscid Euler equations (i.e., Eqs. (\ref{cons_mass}), (\ref{cons_mom}), and (\ref{cons_energy})) coupled to an arbitrary EOS has been performed by Ovsiannikov~\cite{ovsyannikov1} and numerous other authors. As summarized by, for example, Axford~\cite{axford}, the admissible point-groups of Eqs. (\ref{cons_mass}), (\ref{cons_mom}), and (\ref{cons_energy}) conditinally include time translation, space translation, three separate scalings, and a Galilean boost. Projective symmetries are also available under severely restricted geometries and EOS instantiations~\cite{coggeshall1991analytic}. In conjunction with these analyses, Ovsiannikov~\cite{ovsyannikov1} demonstrates that the time translation and one of the scaling symmetries are always present in Eqs. (\ref{cons_mass}), (\ref{cons_mom}), and (\ref{cons_energy}), regardless of geometry and the form of the EOS. 

In addition to these general studies, Ramsey and Baty~\cite{Ramsey_Baty} provide a complementary discussion centered on the conditions under which Eqs. (\ref{cons_mass}), (\ref{cons_mom}), and (\ref{cons_energy}) are invariant under all three possible scaling groups. As discussed extensively by Barenblatt~\cite{barenblatt1,barenblatt2}, scaling phenonena are of particular importance for numerous reasons: including as manifestations of "phenomenon of basic importance," intermediate asymptotic entities, or guides for the construction of scaled experiments in the appropriate contexts. As such, in the sprit of the work of Barenblatt~\cite{barenblatt1,barenblatt2}, Ramsey and Baty~\cite{Ramsey_Baty}, and Albright et al.~\cite{Albright_Ramsey}, the emphasis of this work is restricted from the broader symmetry classes discussed above to invariance under the scaling transformations indicated by Eqs.~(\ref{tnew})-(\ref{Inew}).

Put simply, if we substitute Eqs.~(\ref{tnew})-(\ref{Inew}) into Eqs.~(\ref{mu_1}) and (\ref{mu_2}), invariance demands that the resulting relations are unchanged aside from the indexing from the original variable to the ``new" variable. For example, the invariance condition for the homentropic inviscid Euler EDS is expressed as
\begin{eqnarray}
\mu_1\left(t_{\rm{new}},r_{\rm{new}},\rho_{\rm{new}},u_{\rm{new}},P_{\rm{new}},I_{\rm{new}}\right) &=& \mu_1\left(t,r,\rho,u,P,I\right),\label{mu1_new_to_mu1}\nonumber\\\\
\mu_2\left(t_{\rm{new}},r_{\rm{new}},\rho_{\rm{new}},u_{\rm{new}},P_{\rm{new}},I_{\rm{new}}\right) &=& \mu_2\left(t,r,\rho,u,P,I\right).\label{mu2_new_to_mu2}\nonumber\\
\end{eqnarray}
As originally demonstrated by Sophus Lie~\cite{lie1888theorie1,lie1890theorie2,lie1891vorlesungen,lie1893vorlesungen,lie1893theorie3,lie1896theorie}, this global concept of invariance may be equivalently realized in terms of a local (or infinitesimal) representation in terms of a Lie derivative operation. While potentially not immediately recognizable as advantageous in the context of Eqs.~(\ref{tnew})-(\ref{Inew}), as was the case with the EDS representation discussed in Sec.~\ref{diff_forms}, Lie's formalism is host to distinct advantages. Namely, when generalized to arbitrary transformations, the equations that determine the form of Eqs.~(\ref{tnew})-(\ref{Inew}) leaving Eqs.~(\ref{mu1_new_to_mu1}) and (\ref{mu2_new_to_mu2}) invariant are typically nonlinear, and thus may be difficult if not impossible to solve. The infinitesimal framework reduces all determining equations for $a_1 - a_6$ (or their generalization to arbitrary transformations) to linear equations. 

To construct the infinitesimal analog of Eqs.~(\ref{mu1_new_to_mu1}) and (\ref{mu2_new_to_mu2}), the left-hand sides of these relations are expanded in a Taylor series about the identity element $\epsilon = 0$:
\begin{equation}\label{mu_i_new_taylor_expanded}
\mu_{i,{\rm{new}}} = \mu_i + \epsilon\frac{\partial \mu_i}{\partial\epsilon}\biggr\rvert_{\epsilon=0}+ \frac{1}{2}\epsilon^2\frac{\partial^2\mu_i}{\partial\epsilon^2}\biggr\rvert_{\epsilon=0}+\cdots.
\end{equation}
Using the chain rule, we find the $\epsilon$-derivative may be re-expressed as:
\begin{eqnarray}\label{epsilon_derivative}
\frac{\partial}{\partial\epsilon} &=& \frac{\partial t_{\rm{new}}}{\partial\epsilon}\biggr\rvert_{\epsilon=0}\frac{\partial}{\partial t} +\frac{\partial r_{\rm{new}}}{\partial\epsilon}\biggr\rvert_{\epsilon=0}\frac{\partial}{\partial r}+\frac{\partial\rho_{\rm{new}}}{\partial\epsilon}\biggr\rvert_{\epsilon=0}\frac{\partial}{\partial \rho}\nonumber\\
&&+\frac{\partial u_{\rm{new}}}{\partial\epsilon}\biggr\rvert_{\epsilon=0}\frac{\partial}{\partial u}+\frac{\partial P_{\rm{new}}}{\partial\epsilon}\biggr\rvert_{\epsilon=0}\frac{\partial}{\partial P}+\frac{\partial I_{\rm{new}}}{\partial\epsilon}\biggr\rvert_{\epsilon=0}\frac{\partial}{\partial I}.\nonumber\\
\end{eqnarray}
Applying the appropriate derivatives to Eqs.~(\ref{tnew})-(\ref{Inew}) and substituting them into Eq.~(\ref{epsilon_derivative}) we find
\begin{eqnarray}\label{scaling_group_generator}
\frac{\partial}{\partial\epsilon} = \chi &=& a_1 t \frac{\partial}{\partial t} + a_2 r \frac{\partial}{\partial r} + a_3 \rho \frac{\partial}{\partial \rho} + a_4 u \frac{\partial}{\partial u}+a_5 P \frac{\partial}{\partial P} \nonumber\\
&& + a_6 I \frac{\partial}{\partial I}.
\end{eqnarray}
Using the results from Eq.~(\ref{scaling_group_generator}), Eq.~(\ref{mu_i_new_taylor_expanded}) becomes 
\begin{equation}
\mu_{i,{\rm{new}}} = \mu_i+\epsilon\chi\mu_i+\frac{1}{2}\epsilon^2\chi^2\mu_i +\cdots,
\end{equation}
and using Eq.~(\ref{mu1_new_to_mu1}) and (\ref{mu2_new_to_mu2}) (i.e., the global invariance condition), we find
\begin{equation}\label{final_mu_condition}
\epsilon\chi\mu_i+\frac{1}{2}\epsilon^2\chi^2\mu_i +\cdots =0.
\end{equation}
Therefore, the nontrivial (i.e., $\epsilon \neq 0$) solution of Eq.~(\ref{final_mu_condition}) is
\begin{eqnarray}\label{generator_on_mu_i}
a_1t\frac{\partial\mu_i}{\partial t} + a_2r\frac{\partial\mu_i}{\partial r} + a_3\rho \frac{\partial\mu_i}{\partial \rho} + a_4 u \frac{\partial\mu_i}{\partial u}+a_5 P \frac{\partial\mu_i}{\partial P}\nonumber\\+a_6 I \frac{\partial\mu_i}{\partial I}=0,\;\;\;\;\;\;\;\;\;\; 
\end{eqnarray}
if and only if
\begin{equation}\label{mu_i_zero}
\mu_i =0,
\end{equation}
for all $i$. This infinitesimal invariance condition is entirely equivalent to the global invariance condition given by Eqs.~(\ref{mu1_new_to_mu1}) and (\ref{mu2_new_to_mu2}).

Finally, while Eq.~(\ref{generator_on_mu_i}) represents invariance only of the EDS system, any ancillary conditions appearing in a problem formulation must similarly be invariant under the operation of the Lie derivative or ``group generator" $\chi$.
% 
%%%%%%%%%%%%%%%%%%%%%%%%%%%%%%%%%%%%%%%%%%%
%				          %
%               ANALYSIS                  %
%				          %
%%%%%%%%%%%%%%%%%%%%%%%%%%%%%%%%%%%%%%%%%%%
%
%
\section{Scaling Analysis}\label{scaling_analysis}

Having constructed the scaling group generator $\chi$ given by Eq.~(\ref{scaling_group_generator}), all features of a given problem must be simultaneously invariant under its operation for the entire problem to be invariant under the indicated group of scaling transformations. For a problem featuring a shock wave in a fluid characterized by an isentropic EOS, not only must Eqs.~(\ref{mu_1}) and (\ref{mu_2}) be invariant (as indicated by Eqs.~(\ref{generator_on_mu_i}) and (\ref{mu_i_zero})), but so must:
\begin{itemize}
\item The shock jump conditions given by Eqs.~(\ref{rankine_hugoniot_mass_jump}) and (\ref{rankine_hugoniot_mom_jump}),
\item The conditions on the isentropic EOS given by Eqs.~(\ref{eos}), (\ref{Isolve}), and (\ref{law}),
\item Any other conditions specific to a problem under investigation.
\end{itemize}

If at least one of the scaling constants $a_i$ appearing in Eq.~(\ref{scaling_group_generator}) is revealed to be non-zero as an outcome of the analysis, then the problem is invariant under a scaling transformation.

\subsection{Reduced Inviscid Euler Equations}\label{euler}

In evaluating Eq.~(\ref{generator_on_mu_i}) with Eq.~(\ref{mu_i_zero}), it is necessary to understand the interaction of the group generator $\chi$ with both the exterior derivative $d$ and wedge product $\wedge$ operators. As noted in Sec.~\ref{invariance}, the group generator $\chi$ is actually a Lie derivative, which is itself a generalization in a space of arbitrary dimension of the more familiar directional derivative as appearing in elementary vector calculus. As discussed by Edelen~\cite{edelen2005applied}, Suhubi~\cite{suhubi2013exterior}, and many others, Lie and exterior derivatives commute:
\begin{equation}
\chi d \left(q_i\right) = d \chi\left(q_i\right), \label{commute}
\end{equation}
where $q_i$ retains its previous definition. Moreover, the Lie derivative applied to an arbitrary 2-form obeys the product rule for derivatives:
\begin{equation}
\chi \left(dq_1 \wedge dq_2\right) = d\chi q_1 \wedge q_2 + q_1 \wedge d \chi q_2, \label{product}
\end{equation}
and easily generalizes to forms of arbitrary order. These two important properties are another example of of an advantage of the isovector formalism: group generator operations on differential forms are simple to evaluate.

Using the properties from Eqs.~(\ref{commute}) and~(\ref{product}), Eq.~(\ref{mu_1}), with Eq.~(\ref{mu_i_zero}) for $i = 1$, becomes
\begin{widetext}
\begin{eqnarray}\label{chi_mu1}
\chi\mu_1&=&\left(\chi r\right)d\rho\wedge dr + r d\left(\chi\rho\right)\wedge dr + r d\rho\wedge d\left(\chi r\right)-\left(\chi r\right)u d\rho\wedge dt - r\left(\chi u\right)d\rho\wedge dt- ru d\left(\chi \rho\right)\wedge dt \nonumber\\ 
&&- ru d\rho\wedge d\left(\chi t\right)-\left(\chi \rho\right)rdu\wedge dt -\rho\left(\chi r\right)du\wedge dt-\rho r d\left(\chi u\right)\wedge dt -\rho r du\wedge d\left(\chi t\right)\nonumber\\
&&+n\left(\chi\rho\right)u dt\wedge dr+n\rho \left(\chi u\right)dt\wedge dr+n\rho ud\left(\chi t\right)\wedge dr +n\rho u dt\wedge d\left(\chi r\right) =0.
\end{eqnarray}
\end{widetext}
Applying each derivative of the group generator leads to
\begin{eqnarray}
&&\left(a_2+a_3+a_2\right)r d\rho\wedge dr - \left(a_2+a_4+a_3+a_1\right)ru d\rho\wedge dt \nonumber\\
&&-\left(a_2+a_3+a_4+a_1\right)r\rho du\wedge dt\nonumber\\
&&+ \left(a_3+a_4+a_1+a_2\right)n\rho u dt\wedge dr =0,
\end{eqnarray}
which, using Eq.~(\ref{mu_1}), simplifies to
\begin{eqnarray}
&&\left(a_2+a_3+a_2\right)\left(rud\rho\wedge dt + r\rho du\wedge dt - r\rho u dt\wedge dr\right) \nonumber\\
&&- \left(a_2+a_4+a_3+a_1\right)ru d\rho\wedge dt \nonumber\\
&&-\left(a_2+a_3+a_4+a_1\right)r\rho du\wedge dt\nonumber\\
&&+ \left(a_3+a_4+a_1+a_2\right)n\rho u dt\wedge dr =0.
\end{eqnarray}
Multiplying each term and simplifying, we have
\begin{eqnarray}\label{completed_chi_mu1}
&&\left(a_2+a_3+a_2-a_2-a_4-a_3-a_1 \right)ru d\rho\wedge dt\nonumber\\
&&+ \left(a_2+a_3+a_2-a_3-a_4-a_1-a_2\right)r\rho du\wedge dt\nonumber\\
&&+\left(a_3+a_4+a_1+a_2-a_2-a_3-a_2\right)n\rho u dt\wedge dr =0.\nonumber\\
\end{eqnarray}
For Eq.~(\ref{completed_chi_mu1}) to be nontrivially satisfied, the coefficient of each unique 2-form appearing within it must be zero. This procedure yields three redundant conditions:
\begin{equation}\label{a_4_constraining_equation}
a_4=a_2-a_1,
\end{equation}
which encapsulates the dimensionally correct statement that the fluid velocity $u$ scales as $\frac{r}{t}$.

We can perform a similar analysis on Eq.~(\ref{mu_2}):
\begin{widetext}
\begin{eqnarray}
\chi\mu_2&=&\left(\chi \rho^2\right)du\wedge dr + \rho^2 d\left(\chi u\right)\wedge dr + \rho^2 du\wedge d\left(\chi r\right)-\left( \chi \rho^2 \right)u du\wedge dt - \rho^2 \left( \chi u\right)du\wedge dt\nonumber\\
&&-\rho^2u d\left( \chi u\right)\wedge dt -\rho^2u du\wedge d\left( \chi t\right)-\left[ \chi K_S(\rho)\right] d\rho\wedge dt -K_S d\rho \wedge d\left( \chi t\right)=0.
\end{eqnarray}
\end{widetext}
Simplifying, as we have done above, we have
\begin{eqnarray}\label{completed_chi_mu2}
&&\left(-a_2+a_4+a_1\right)\rho^2 du \wedge dr \nonumber\\
&&+ \left[K_S \left(-a_3-2a_4\right) +a_3 \rho \frac{d K_S}{d\rho}\right] d\rho \wedge dt = 0.
\end{eqnarray}
Again, for Eq.~(\ref{completed_chi_mu2}) to be nontrivially satisfied, the coefficient of each unique 2-form appearing within it must be zero. The coefficient of the $du\wedge dr$ term is identically zero by Eq. (\ref{a_4_constraining_equation}), leaving for the $d\rho\wedge dt$ term:
\begin{eqnarray}\label{final_bulk_ode}
a_3 \rho \frac{d K_S}{d\rho} - \left(a_3+2a_2-2a_1\right)K_S = 0,
\end{eqnarray}
which is an ordinary differential equation (ODE) for $K_S$ that has multiple solutions. 

Given Eq.~(\ref{a_4_constraining_equation}), Eq.~(\ref{generator_on_mu_i}) becomes
\begin{eqnarray}\label{final_euler_group_generator} 
\chi &=& a_1 t \frac{\partial}{\partial t} + a_2 r \frac{\partial}{\partial r} + a_3 \rho \frac{\partial}{\partial \rho} + \left(a_2-a_1\right) u \frac{\partial}{\partial u}+a_5 P \frac{\partial}{\partial P}\nonumber\\
&&+a_6 I \frac{\partial}{\partial I}.
\end{eqnarray}
This result is the most general group generator associated with invariance of the homentropic inviscid Euler equations. However, this generator may assume different forms depending on the choice of $K_S$.

%%%%%%%%%%%%%%%%%%%%%%%%%%%%%%%%%%%%%%%%%%%%%%%%%%%%
%					           %
%   Invariance of Rankine-Hugoniot jump conditions %
%						   %
%%%%%%%%%%%%%%%%%%%%%%%%%%%%%%%%%%%%%%%%%%%%%%%%%%%%
\subsection{Rankine-Hugoniot jump conditions}\label{jump}

We now perform the same analysis as we did in Sec.~\ref{euler} for Eqs.~(\ref{rankine_hugoniot_mass_jump}) and (\ref{rankine_hugoniot_mom_jump}) by applying Eq.~(\ref{final_euler_group_generator})\footnote{For completeness, we use Eq. (\ref{final_euler_group_generator}) in order to find the constraints on all variables.}. First we rewrite Eqs.~(\ref{rankine_hugoniot_mass_jump}) and (\ref{rankine_hugoniot_mom_jump}) as follows:
\begin{eqnarray}
\mu_3 &\equiv& \left(u_1-D\right)\rho_1 - \left(u_2-D\right)\rho_2,\label{mu3}\\
\mu_4 &\equiv& P_1- P_2 +\rho_1\left(u_1-D\right)u_1 - \rho_2\left(u_2-D\right)u_2.\label{mu4}
\end{eqnarray}
Applying Eq.~(\ref{final_euler_group_generator}) to Eq.~(\ref{mu3}) we have
\begin{eqnarray}
\chi\mu_3&=&\left(a_3+a_2-a_1\right)\left(\rho_1u_1-\rho_2u_2\right)+\left(a_1t\frac{dD}{dt}+a_3D\right)\rho_2\nonumber\\
&&-\left(a_1t\frac{dD}{dt}+a_3D\right)\rho_1=0.
\end{eqnarray}

Using Eq.~(\ref{mu3}) we have:
\begin{eqnarray}
&&\left[a_1t\frac{dD}{dt}+\left(a_3-a_3-a_2+a_1\right)D\right]\rho_2\nonumber\\
&&-\left[a_1t\frac{dD}{dt}+\left(a_3-a_3-a_2+a_1\right)D\right]\rho_1=0.
\end{eqnarray}

To satisfy this invariance condition, we must solve the following differential equation:
\begin{eqnarray}\label{diff_eq_d}
a_1t\frac{dD}{dt}+\left(a_1-a_2\right)D=0,
\end{eqnarray}
whose nontrivial (i.e., $a_1 \neq 0$) solution is
\begin{eqnarray}\label{D_solution}
D\left(t\right)=D_0t^\sigma,
\end{eqnarray}
where $\sigma=\frac{a_2-a_1}{a_1}$, and $D_0$ is an arbitrary integration constant.

We now apply Eq.~(\ref{final_euler_group_generator}) to Eq.~(\ref{mu4}) which results in
\begin{widetext}
\begin{eqnarray}
\chi\mu_4&=&a_5\left(P_1-P_2\right)+\left(a_3+2a_2-2a_1\right)\rho_1u_1^2-\left[\left(a_3+a_2-a_1\right)D+a_1 t \frac{dD}{dt}\right]\nonumber\\
&&+\left(a_3+2a_2-2a_1\right)\rho_1u_1^2-\left[\left(a_3+a_2-a_1\right)D+a_1 t \frac{dD}{dt}\right].
\end{eqnarray}
\end{widetext}
Using Eq.~(\ref{mu4}) and collecting like terms, we have the following determining equations:
\begin{eqnarray}
a_5-a_3-2a_2+2a_1&=&0,\\
a_1t\frac{dD}{dt}+\left(a_1-a_2\right)D&=&0.\label{d_eq_1}
\end{eqnarray}
Since Eq.~(\ref{d_eq_1}) has been satisfied by Eq.~(\ref{D_solution}) we find that 
\begin{equation}\label{a5_condition}
a_5=a_3+2a_2-2a_1,
\end{equation}
which further reduces Eq.~(\ref{final_euler_group_generator}) to
\begin{eqnarray}\label{final_group_generator} 
\chi &=& a_1 t \frac{\partial}{\partial t} + a_2 r \frac{\partial}{\partial r} + a_3 \rho \frac{\partial}{\partial \rho} + \left(a_2-a_1\right) u \frac{\partial}{\partial u}\nonumber\\
&&+\left(a_3+2a_2-2a_1\right)P \frac{\partial}{\partial P}+a_6 I \frac{\partial}{\partial I}.
\end{eqnarray}
Like Eq. (\ref{a_4_constraining_equation}), Eq. (\ref{a5_condition}) encapsulates the dimensionally correct statement that pressure $P$ scales like the density $\rho$ times a specific energy, which has units of $r^2/t^2$.
%
%
%%%%%%%%%%%%%%%%%%%%%%%%%%%%%%%%%%%%%%%%%%
%					 %
%   Invariance of Auxiliary Conditions   %
%					 %
%%%%%%%%%%%%%%%%%%%%%%%%%%%%%%%%%%%%%%%%%%
\subsection{Thermodynamic Constraints}

In addition to the dynamical equations, the thermodynamic constraints encoded in the homentropic inviscid Euler system must likewise be invariant, as has already been partially established in the construction of Eq.~(\ref{final_bulk_ode}). Further constraints include the fundamental thermodynamic relation given by Eq.~(\ref{law}) with $dS = 0$ and the isentropic definitions of $P$ and $I$ given by Eqs.~(\ref{eos}) and (\ref{Isolve}):
\begin{eqnarray}
\mu_5 &=& dI - \frac{P}{\rho^2}d\rho, \label{mu5}\\
\mu_6 &=& P - f\left(\rho\right), \label{mu6}\\
\mu_7 &=& I - g\left(\rho\right) - I_0. \label{mu7}
\end{eqnarray}
We first conduct our symmetry analysis on Eq.~(\ref{mu5}) with the understanding that the fundamental thermodynamic relation is in fact a 1-form, and the total derivatives appearing within it may be regarded as exterior derivatives. As such, invariance of this relation under $\chi$ demands
\begin{equation}\label{gibbs_inv}
\chi \mu_5 = a_6 dI - \left(2a_2 - 2a_1\right) \frac{P}{\rho^2} d\rho = 0,
\end{equation}
when $\mu_5 = 0$, which yields
\begin{equation}\label{a6_constrained}
a_6 = 2a_2 - 2a_1.
\end{equation}
Again, this reinfores our prior understanding that specific energy $I$ scales in units of $r^2/t^2$.

As a result, Eq.~(\ref{final_group_generator}) can now be written entirely in terms of the scaling constants $a_1$, $a_2$, and $a_3$:
\begin{eqnarray}\label{final_final_group_generator} 
\chi &=& a_1 t \frac{\partial}{\partial t} + a_2 r \frac{\partial}{\partial r} + a_3 \rho \frac{\partial}{\partial \rho} + \left(a_2-a_1\right) u \frac{\partial}{\partial u}\nonumber\\
&&+\left(a_3+2a_2-2a_1\right)P \frac{\partial}{\partial P}+\left(2a_2-2a_1\right) I \frac{\partial}{\partial I}.\nonumber\\
\end{eqnarray}

With Eq.~(\ref{final_final_group_generator}), invariance of the isentropic definition of the pressure $P$ (i.e. Eq.~(\ref{mu6})) demands
\begin{equation}\label{isen_P_inv}
\chi \left[P - f\left(\rho\right)\right] = \left(a_3+2a_2-2a_1\right)P - a_3 f'\left(\rho\right) = 0,
\end{equation}
when $\mu_6 = 0$, which yields an ODE that $f$ must satisfy:
\begin{equation}\label{f_constrained}
\left(a_3+2a_2-2a_1\right)f\left(\rho\right) - a_3 \rho f'\left(\rho\right) = 0.
\end{equation}

Likewise, with Eq.~(\ref{final_final_group_generator}), invariance of the isentropic definition of the SIE $I$ (i.e. Eq.~(\ref{mu7})) demands
\begin{equation}\label{isen_I_inv}
\chi \left[I - g\left(\rho\right) - I_0\right] = \left(2a_2-2a_1\right)I - a_3 g'\left(\rho\right) = 0,
\end{equation}
when $\mu_7 = 0$, which yields an ODE that $g$ must satisfy:
\begin{equation}\label{g_constrained}
\left(2a_2-2a_1\right)\left[g\left(\rho\right)+I_0\right] - a_3 \rho g'\left(\rho\right) = 0.
\end{equation}

The EOS functions $f$ and $g$ are connected through satisfaction of the (isentropic) fundamental thermodynamic relation represented by Eq.~(\ref{gdef}), and must also be consistent with an associated adiabatic bulk modulus $K_S$ calculated via Eq.~(\ref{abm_isen}). As a result, the solutions of the ODEs given by Eqs.~(\ref{final_bulk_ode}), (\ref{f_constrained}), and (\ref{g_constrained}) must be mutually consistent so as to enable the presence of various scaling symmetries. The possible solutions of these three ODEs fall under four cases.

%%%%%%%%%%%%%%%%%%%%%%%%%%%%%%%%%%%%%%%%%%%%%%%%%%
%						 %
%		     Case I 			 % 
%						 %
%%%%%%%%%%%%%%%%%%%%%%%%%%%%%%%%%%%%%%%%%%%%%%%%%%
\subsubsection{Case I: $a_1\neq a_2\neq a_3 \neq 0$}\label{case_1}
The solution to Eq.~(\ref{final_bulk_ode}) is 
\begin{equation}\label{bulk_mod_solution}
K_S\left(\rho\right) = A_1\rho^\psi,
\end{equation}
where $A_1$ is an arbitrary integration constant, and $\psi\equiv\frac{a_3+2a_2-2a_1}{a_3}$. With Eqs.~(\ref{abm_isen}) and (\ref{Isolve}), the associated EOS for $P$ and $I$ is given by
\begin{eqnarray}
P &=& \frac{A_1}{\psi}\rho^\psi + P_0 , \label{P_case1} \\
I &=& \frac{A_1}{\psi \left(\psi-1\right)} \rho^{\psi - 1} -\frac{P_0}{\rho} + I_0 , \label{I_case1}
\end{eqnarray}
where $P_0$ is an arbitrary integration constant. Inserting these results into Eqs.~(\ref{f_constrained}) and (\ref{g_constrained}) results in the requirements $P_0 = I_0 = 0$ for these constraints to be satisfied; $A_1$ and $\psi$ are otherwise unconstrained (aside from assuming values that yield positive $P$ and $I$, on the grounds of physical realism).

Moreover, with $P_0 = I_0 = 0$, Eqs.~(\ref{P_case1}) and (\ref{I_case1}) may be combined to yield an EOS of the form $P = P\left(\rho,I\right)$:
\begin{equation}\label{P-rho-I_case1}
P = \left(\psi -1\right) \rho I,
\end{equation}
which is of the ideal gas type. This case is thus associated with the wide body of existing literature associated with scaling solutions in the context of an ideal gas EOS~\cite{sedov,zeldovich_raizer}.

In this case, Eq.~(\ref{final_final_group_generator}) is as indicated.
%%%%%%%%%%%%%%%%%%%%%%%%%%%%%%%%%%%%%%%%%%%%%%%%%%
%						 %
%		     Case II 			 % 
%		 				 %
%%%%%%%%%%%%%%%%%%%%%%%%%%%%%%%%%%%%%%%%%%%%%%%%%%
\subsubsection{Case II: $a_1=a_2$ and $a_3\neq 0$}\label{case_2}
We can rewrite Eq.~(\ref{final_bulk_ode}) as follows
\begin{eqnarray}
a_3\rho \frac{\partial K_S}{\partial \rho}-a_3 K_S &=&0,
\end{eqnarray}
whose solution is
\begin{equation}\label{case2_abm}
K_S\left(\rho\right) = A_2\rho,
\end{equation}
where $A_2$ is an arbitrary integration constant. With Eqs.~(\ref{abm_isen}) and (\ref{Isolve}), the associated EOS for $P$ and $I$ is given by
\begin{eqnarray}
P &=& A_2 \rho + P_0, \label{P_case2} \\
I &=& A_2 \ln(\rho) - \frac{P_0}{\rho} + I_0, \label{I_case2}
\end{eqnarray}
where $P_0$ is an arbitrary integration constant. Inserting these results into Eqs.~(\ref{f_constrained}) and (\ref{g_constrained}) results in the requirements $A_2 = 0$ and $P_0 = 0$ for these constraints to be satisfied; $I_0$ is otherwise unconstrained. This case is thus trivial; as it features $P = 0$ for any $\rho$ and $I$.

In this case, Eq.~(\ref{final_final_group_generator}) reduces to
\begin{equation}\label{final_case_II_group_generator}
\chi = a_1 t \frac{\partial}{\partial t} + a_1 r \frac{\partial}{\partial r} + a_3 \rho \frac{\partial}{\partial \rho} +a_3 P \frac{\partial}{\partial P}.
\end{equation}
%%%%%%%%%%%%%%%%%%%%%%%%%%%%%%%%%%%%%%%%%%%%%%%%%%
%						 %
%		     Case III 			 % 
%		 				 %
%%%%%%%%%%%%%%%%%%%%%%%%%%%%%%%%%%%%%%%%%%%%%%%%%%
\subsubsection{Case III: $a_1\neq a_2$ and $a_3= 0$}\label{case_3}
We can rewrite Eq.~(\ref{final_bulk_ode}) as follows
\begin{equation}
\left(a_1-a_2\right)K_S =0, \label{case3_bulk_ode}
\end{equation}
whose solution is $K_S(\rho)=0$. With Eqs.~(\ref{abm_isen}) and (\ref{Isolve}), the associated EOS for $P$ and $I$ is given by
\begin{eqnarray}
P &=& P_0, \label{P_case3} \\
I &=& -\frac{P_0}{\rho} + I_0, \label{I_case3}
\end{eqnarray}
where $P_0$ is an arbitrary integration constant. Inserting these results into Eqs.~(\ref{f_constrained}) and (\ref{g_constrained}) results in the requirements $P_0 = I_0 = 0$ for these constraints to be satisfied. Like the previous case, this case is trivial as it features $P = I = 0$ for any $\rho$.

In this case, Eq.~(\ref{final_final_group_generator}) reduces to
\begin{eqnarray}\label{final_case_III_group_generator}
\chi &=& a_1 t \frac{\partial}{\partial t} + a_2 r \frac{\partial}{\partial r} + \left(a_2-a_1\right) u \frac{\partial}{\partial u}+\left(2a_2-2a_1\right)P \frac{\partial}{\partial P}\nonumber\\
&&+\left(2a_2-2a_1\right) I \frac{\partial}{\partial I}.
\end{eqnarray}
%%%%%%%%%%%%%%%%%%%%%%%%%%%%%%%%%%%%%%%%%%%%%%%%%%
%						 %
%		     Case IV 			 % 
%						 %
%%%%%%%%%%%%%%%%%%%%%%%%%%%%%%%%%%%%%%%%%%%%%%%%%%
\subsubsection{Case IV: $a_1=a_2$ and $a_3=0$}\label{case_4}
We can see that Eq.~(\ref{final_bulk_ode}) is solved identically, allowing the adiabatic bulk modulus and associated EOS to be unconstrained aside from the thermodynamic requirements given by Eqs.~(\ref{abm_isen}) and (\ref{Isolve}). This case is thus a direct manifestation of the ``universal" scaling symmetry as discussed by Ovsiannkov~\cite{ovsyannikov1} and Ramsey and Baty~\cite{Ramsey_Baty}; physically, this symmetry is associated with purely kinematic (i.e., $r-t$ only) scaling of the space and time variables, and not any dependent variables. As such, on physical grounds the EOS is necessarily unconstrained.

In this case, Eq.~(\ref{final_final_group_generator}) becomes
\begin{eqnarray}\label{case_4_group_generator}
\chi = a_1 t \frac{\partial}{\partial t} + a_1 r \frac{\partial}{\partial r}.
\end{eqnarray}

To summarize, the four possible EOS cases appearing in Secs.~\ref{case_1}-\ref{case_4} are reported in Table~\ref{table1}.
\renewcommand{\arraystretch}{1.5}
\begin{table}[t]
\begin{tabular}{c c c c c} 
 \hline
  & $K_S(\rho)$ & $P(\rho)$ & $I(\rho)$ & $D(t)$ \\ [1.0ex] 
 \hline\hline
 $a_1\neq a_2\neq a_3 \neq 0$ & $A_1\rho^\psi$ & $\frac{A_1}{\psi}\rho^\psi$ &    $\frac{A_1}{\psi\left(\psi-1\right)}\rho^{\psi-1}+I_0$ & $D_0t^\sigma$  \\ 
  $a_1=a_2$, $a_3\neq 0$ & $0$ & $0$ & $I_0$ & $D_0$\\
 $a_1\neq a_2$, $a_3 =0$ & $0$ & $0$ & $0$ & $D_0t^\sigma$\\
 $a_1=a_2$, $a_3 = 0$& Arbitrary & Arbitrary & Arbitrary & $D_0$\\
 \hline
\end{tabular}
\caption{Summary of possible scale-invariant isentropic EOS classes and associated shock trajectories. Both $f(\rho)$ and $g(\rho)$ are replaced with their physical variables $P(\rho)$ and $I(\rho)$, respectively.}
\label{table1}
\end{table}
%%%%%%%%%%%%%%%%%%%%%%%%%%%%%%%%%%%%%%%%%%
%					 %
%	    Similarity Variables	 %
%					 %
%%%%%%%%%%%%%%%%%%%%%%%%%%%%%%%%%%%%%%%%%%
%\vspace{-0.6cm}
\subsection{Similarity Variables}\label{similarity_variables}
In addition to providing a convenient means for interrogating the group invariance properties of various algebraic structures, infinitesimal group generators can also be used to construct changes of coordinates in terms of which invariant differential equations may be reduced to simpler structures (e.g., from PDEs to ODEs). In particular, for an arbitrary function $H$ of all independent and dependent variables spanning a problem formulation (in the current case, comprised of $r, t, \rho, u, P$, and $I$), the PDE condition
\begin{equation}
\label{inv_func}
\chi H\left(r,t,\rho,u,P,I\right) = 0 ,
\end{equation}
indicates that $H$ is invariant under the group of transformations generated by $\chi$. The PDE given by Eq.~(\ref{inv_func}) may be solved using the Method of Characteristics. The arbitrary constants of integration arising from this solution are invariant by construction under the action of $\chi$, and thus may be used to develop similarity variables in terms of which the original system of PDEs may be reformulated and simplified.

The characteristic equations associated with Eq.~(\ref{final_final_group_generator}) are given by
\begin{widetext}
\begin{eqnarray}
\frac{dt}{a_1t}=\frac{dr}{a_2r}=\frac{d\rho}{a_3\rho}=\frac{du}{\left(a_2-a_1\right)u}=\frac{dP}{\left(a_3+2a_2-2a_1\right)P}=\frac{dI}{\left(2a_2-2a_1\right) I}=\frac{dH}{0},
\end{eqnarray}
\end{widetext}
or, $H = \rm{const.}$ along the characteristic curves defined by

\begin{eqnarray}
\frac{dt}{a_1t}&=&\frac{dr}{a_2r} \label{dt_dr_sim},\\
\frac{d\rho}{a_3\rho}&=&\frac{dr}{a_2r} \label{dr_drho_sim},\\
\frac{du}{\left(a_2-a_1\right)u}&=&\frac{dr}{a_2r} \label{dr_du_sim},\\
\frac{dP}{\left(a_3+2a_2-2a_1\right)P}&=&\frac{dr}{a_2r} \label{dr_dP_sim},\\
\frac{dI}{\left(2a_2-2a_1\right)I}&=&\frac{dr}{a_2r} \label{dr_dI_sim}.
\end{eqnarray}

The solutions of Eqs.~(\ref{dt_dr_sim})-(\ref{dr_dI_sim}) are 

\begin{eqnarray}
\xi &=& \frac{r}{t^\alpha}, \label{xi_def}\\
\rho &=& r^\zeta w\left(\xi\right), \label{rho_sim_sol}
\end{eqnarray}
\begin{eqnarray}
P&=&r^\lambda m\left(\xi\right), \label{P_sim_sol}\\
u&=&r^\beta j\left(\xi\right), \label{u_sim_sol}\\
I&=&r^\tau h\left(\xi\right), \label{I_sim_sol}
\end{eqnarray}
where $\xi$, $w$, $m$, $j$, and $h$ are the constants of integration (i.e., the invariants of the group) that may be interpreted as a change of variables, and $\alpha\equiv\frac{a_2}{a_1}$, $\zeta\equiv\frac{a_3}{a_2}$, $\tau\equiv\frac{2a_2-2a_1}{a_2}$, $\lambda\equiv\frac{a_3+2a_2-2a_1}{a_2}$, and $\beta\equiv \frac{a_2-a_1}{a_2}$ (as summarized in Table~\ref{table2}, along with other constants appearing elsewhere).
\renewcommand{\arraystretch}{1.5}
\begin{table}[t]
\centering
\begin{tabular}{||c c||} 
 \hline
  Power Law Variable & Scaling Constants \\ [1.0ex] 
 \hline\hline
 $\sigma$& $\frac{a_2-a_1}{a_1}$\\
 $\psi$& $\frac{a_3+2a_2-2a_1}{a_3}$ \\
 $\alpha $ & $\frac{a_2}{a_1}$ \\ 
 $\zeta$& $\frac{a_3}{a_2} $ \\
 $\lambda$& $\frac{a_3+2a_2-2a_1}{a_2} $\\
 $\beta$ & $\frac{a_2-a_1}{a_2}$ \\
 $\tau$& $ \frac{2a_2-2a_1}{a_2}$ \\
 \hline
\end{tabular}
\caption{Summary of variable definitions involving scaling constants.}
\label{table2}
\end{table}
As detailed in Appendix~\ref{sim_sub_euler}, we now substitute Eqs.~(\ref{rho_sim_sol}) and (\ref{u_sim_sol}) into Eqs.~(\ref{simp_eos_cons_mass}) and (\ref{simp_eos_cons_mom}) and find the following coupled, reduced system of equations:
\begin{eqnarray}
w'&=&\frac{\xi w j'+\left(\zeta+\beta+n\right)jw}{\xi\left(\alpha\xi^{\frac{1}{\alpha}}-j\right)},\label{new_rho_w}\\
j'&=&\frac{r^{-\zeta-2\beta}K_S\left[r^\zeta w\right]\left(\zeta w +\xi w' \right)+\beta w^2j^2 }{\xi w^2 \left(\alpha\xi^{\frac{1}{\alpha}}-j\right) }\label{new_u_j},
\end{eqnarray}
where the primes indicate ordinary derivatives with respect to the new independent variable $\xi$, and $K_S$ is a function of the indicated argument in square brackets. Equation~(\ref{new_rho_w}) is an ODE in $w(\xi)$ and $j(\xi)$. As written, Eq.~(\ref{new_u_j}) is not, but it reduces further for each of the four cases outlined in Secs.~\ref{case_1}-\ref{case_4}. 

%
%
%%%%%%%%%%%%%%%%%%%%%%%%%%%%%%%%%%%%%%%%%%%%%%%%%%
%						 %
%		     Case I 			 % 
%						 %
%%%%%%%%%%%%%%%%%%%%%%%%%%%%%%%%%%%%%%%%%%%%%%%%%%
%HIGHLIGHT HERE
\subsubsection{Case I: $a_1\neq a_2\neq a_3 \neq 0$}\label{case_1_reduction}

In this case, $\zeta$, $\beta$, and $\alpha$ are as previously indicated, $K_S$ is given by Eq.~(\ref{bulk_mod_solution}), and Eq.~(\ref{new_rho_w}) is as indicated. Equation~(\ref{new_u_j}) becomes
\begin{eqnarray}\label{ode_j_case1}
j'&=& \frac{A_1 w^\psi\left(\zeta w +\xi w' \right)+\beta w^2j^2 }{\xi w^2 \left(\alpha\xi^{\frac{1}{\alpha}}-j\right) }.
\end{eqnarray}

As noted in Sec.~\ref{case_1}, this case is consistent with an ideal gas EOS. As noted in the wide body of existing literature for this case~\cite{sedov,zeldovich_raizer}, Eqs.~(\ref{new_rho_w}) and (\ref{ode_j_case1}) may be further reduced. In particular, with the change of variables
\begin{eqnarray}
j\left(\xi\right) &=& \xi^{\frac{a_1}{a_2}} J\left(\xi\right), \label{J_def} \\
w\left(\xi\right) &=& \xi^{\frac{a_1 a_3}{\left(a_2-a_1\right) a_2}} W\left(\xi\right), \label{W_def}
\end{eqnarray}
Equations~(\ref{new_rho_w}) and (\ref{ode_j_case1}) may be solved as an algebraic system for the derivatives $W'\left(\xi\right)$ and $J'\left(\xi\right)$ to yield
\begin{eqnarray}
\xi W'&=& \frac{\Delta_1}{\Delta}, \label{Wprime_sol} \\
\xi J'&=& \frac{\Delta_2}{\Delta}, \label{Jprime_sol}
\end{eqnarray}
where
\begin{widetext}
\begin{eqnarray}
\Delta_1 &=& a_1^2\left(a_3+na_2\right)J^2W^{\frac{2a_1}{a_3}+1}-a_2\left(a_2-a_1+a_3+na_2\right)JW^{\frac{2a_1}{a_3}+1}-A_1a_1^2a_3W^{\frac{2a_1}{a_3}+1}\nonumber\\
&&+\frac{a_3\left[\left(a_2-a_1J\right)^2W^{\frac{2a_1}{a_3}}-A_1a_1^2W^{\frac{2a_2}{a_3}}\right]}{a_1-a_2}W\label{delta1}\\
\Delta_2 &=&a_1\left(a_1-a_2\right)J^2\left(a_1J-a_2\right)W^{\frac{2a_1}{a_3}}-A_1a_1\left\{a_1\left[a_1-a_2\left(1+n\right)\right]\right\}JW^{\frac{2a_2}{a_3}}+A_1a_1a_2a_3W^{\frac{2a_2}{a_3}}\nonumber\\
&&+a_1\left[\left(a_2-a_1J\right)^2W^{\frac{2a_1}{a_3}}-A_1a_1^2W^{\frac{2a_2}{a_3}}\right]J,\label{delta2}\\
\Delta &=& a_2\left(a_2- a_1 J\right)^2 W^{\frac{2a_1}{a_2}}-A_1a_{1}^{2} a_2 W^{\frac{2a_2}{a_3}}.\label{delta}
\end{eqnarray}
\end{widetext}
Equations~(\ref{Wprime_sol}) and (\ref{Jprime_sol}) are an autonomous system of first order ODEs, and may thus be rewritten as a single ODE
\begin{equation} \label{dWdJ_final}
\frac{dW}{dJ} = \frac{\Delta_1}{\Delta_2},
\end{equation}
supplemented by a quadrature
\begin{equation} \label{dxidJ_final}
\frac{1}{\xi}\frac{d\xi}{dJ} = \frac{\Delta}{\Delta_2},
\end{equation}
which may be evaluated subsequently once $W\left(J\right)$ has been determined from Eq.~(\ref{dWdJ_final}).

%
%%%%%%%%%%%%%%%%%%%%%%%%%%%%%%%%%%%%%%%%%%%%%%%%%%
%						 %
%		     Cases II and III		 % 
%						 %
%%%%%%%%%%%%%%%%%%%%%%%%%%%%%%%%%%%%%%%%%%%%%%%%%%
%HIGHLIGHT HERE
\subsubsection{Cases II and III: $a_1=a_2$ and $a_3\neq 0$ or $a_1\neq a_2$ and $a_3= 0$}\label{case_2_reduction}

In either of these cases, $P = 0$ as discussed in Secs.~\ref{case_2} and \ref{case_3}. Physically, in this scenario, the absence of driving pressure constrains each fluid particle to move with its initial velocity. In the context of Eqs.~(\ref{simp_eos_cons_mass}) and (\ref{simp_eos_cons_mom}), this behavior manifests through the momentum conservation collapsing to the inviscid Burgers' equation, which in general may be solved in isolation for the velocity field; the associated density field may then be constructed through sequential solution of the mass conservation relation, with the velocity solution as input. 

If in this zero-pressure EOS case it is further prescribed that the associated solutions possess either the Case II or III scaling symmetries, the mass density and velocity fields must then satisfy Eqs.~(\ref{new_rho_w}) and (\ref{new_u_j}) with $K_S = 0$ and $\zeta=0$. In this case Eqs.~(\ref{new_rho_w}) and (\ref{new_u_j}) become
\begin{eqnarray}
w'&=&\frac{\xi w j'+\left(\beta+n\right)jw}{\xi\left(\alpha\xi^{\frac{1}{\alpha}}-j\right)},\label{ode_w_case23}\\
j'&=&\frac{\beta w^2j^2 }{\xi w^2 \left(\alpha\xi^{\frac{1}{\alpha}}-j\right) }\label{ode_j_case23},
\end{eqnarray}
As expected given the aforementioned physical arguments, Eq.~(\ref{ode_j_case23}) may be solved independetly for $j$ (and thus the velocity field), and the resulting solution used to sequentally solve Eq.~(\ref{ode_w_case23}). 

Following this procedure for Case~II (featuring $\alpha = 1$ and $\beta = 0$), Eqs.~(\ref{ode_w_case23}) and (\ref{ode_j_case23}) may be further reduced, and solved exactly to yield
\begin{eqnarray}
w\left(\xi\right)&=&w_0\left(\xi-j_0\right)^{n}\xi^{-n}, \label{w_solution_case_2}\\
j &=& j_0, \label{j_solution_case_2}
\end{eqnarray}
where $j_0$ and $w_0$ are arbitrary constants of integration. 

For Case~III, $\alpha$ and $\beta$ are unconstrained, and Eqs.~(\ref{ode_w_case23}) and (\ref{ode_j_case23}) as written have no known closed-form solution.
%
%%%%%%%%%%%%%%%%%%%%%%%%%%%%%%%%%%%%%%%%%%%%%%%%%%
%						 %
%		     Case IV 			 % 
%						 %
%%%%%%%%%%%%%%%%%%%%%%%%%%%%%%%%%%%%%%%%%%%%%%%%%%

\subsubsection{Case IV: $a_1=a_2$ and $a_3=0$}\label{case_4_reduction}
In this case, $\zeta = 0$, $\beta = 0$, and $\alpha = 1$, $K_S$ is arbitrary, and Eq.~(\ref{new_rho_w}) becomes
\begin{equation}\label{ode_w_case4}
w'=\frac{\xi w j'+njw}{\xi\left(\xi-j\right)}.
\end{equation} 
Equation~(\ref{new_u_j}) becomes
\begin{equation}\label{ode_j_case4}
j'=\frac{K_S\left[ w\right] w' }{ w^2 \left(\xi-j\right) }.
\end{equation}
Equations (\ref{ode_w_case4}) and (\ref{ode_j_case4}) are ODEs in $w(\xi)$ and $j(\xi)$, which may be further reduced under perscription of the adiabatic bulk modulus $K_S[w]$. 
%
%
%%%%%%%%%%%%%%%%%%%%%%%%%%%%%%%%%%%%%%%%%%
%					 %
%	      EXAMPLES		 	 %
%					 %
%%%%%%%%%%%%%%%%%%%%%%%%%%%%%%%%%%%%%%%%%%
\section{Example Solutions}\label{examples}
\begin{figure*}[]
	\begin{subfigure}[t]{0.4\textwidth}
	%\centering
        \includegraphics[scale=0.46]{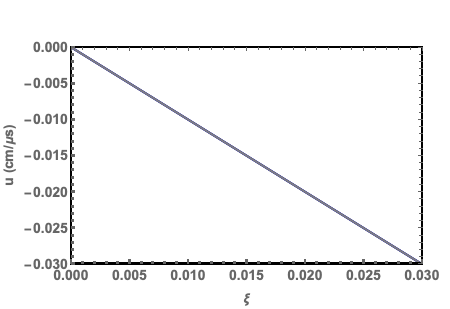}
		\label{smooth_speed_plot}
    \end{subfigure}
        \begin{subfigure}[t]{0.4\textwidth}
 %       \centering
        \includegraphics[scale=0.43]{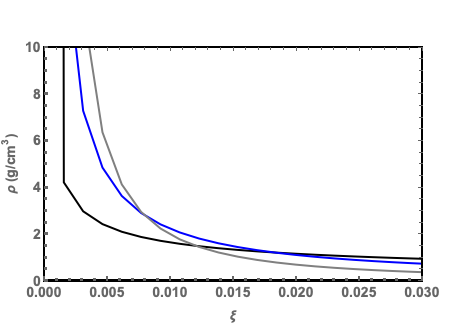}
		\label{smooth_density_plot}
    \end{subfigure}   
    \begin{subfigure}[t]{0.4\textwidth}
 %   \centering
        \includegraphics[scale=0.45]{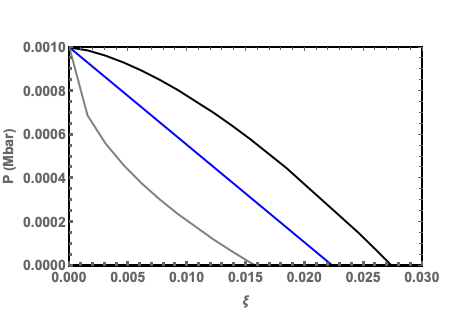}
        \label{smooth_pressure_plot}
    \end{subfigure}
    \begin{subfigure}[t]{0.4\textwidth}
 %   \centering
        \includegraphics[scale=0.43]{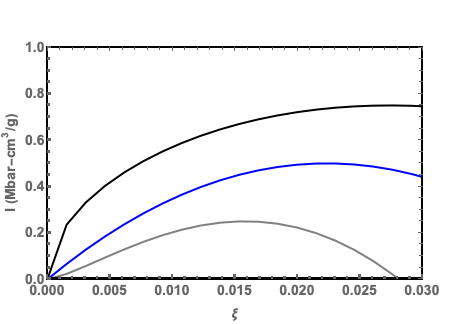}
        \label{smooth_sie_plot}
    \end{subfigure}
\caption{Flow field for example shock-free, homentropic solution with modified Tait EOS. Black, blue, and gray lines denote $n = 0, 1,$ and $2$, respectively. The figures show the speed ($u$), density ($\rho$), pressure ($P$), and SIE ($I$) as a function of similarity variable ($\xi = r/t$).}
\label{smooth_plots}
\end{figure*}

The utility of the the preceding calcluations is demonstrated through their application in constructing exact or semi-analytical solutions of the underlying mathematical model. We provide two examples: (1) a shock-free solution featuring a variety of prescribed properties, and (2) the classical Noh problem, featuring a stagnation shock propagating outward from a rigid wall into a gas infalling with constant velocity.

%%%%%%%%%%%%%%%%%%%%%%%%%%%%%%%%%%%%%%%%%%
%					 %
%       Smooth/Shock-Free Solution 	 %
%					 %
%%%%%%%%%%%%%%%%%%%%%%%%%%%%%%%%%%%%%%%%%%
%HIGHLIGHT HERE
\subsection{Shock-Free Solution}\label{smooth}

For an arbitrary isentropic EOS, scale-invariant compressible flow solutions must originate from Eqs.~(\ref{ode_w_case4}) and (\ref{ode_j_case4}). By construction any such solution of these equations is invariant only under the universal scaling group generated by Eq.~(\ref{case_4_group_generator}), as discussed in Sec.~\ref{case_4}. 

To construct a purely homentropic (i.e., shock-free or smooth) solution of Eqs.~(\ref{ode_w_case4}) and (\ref{ode_j_case4}), consider a flow featuring 
\begin{equation}\label{jnegxi}
j\left(\xi\right) = -\xi .
\end{equation}
In the physical variables associated with Case~IV, this assumption corresponds to 
\begin{equation}\label{speed_smooth_condition}
u\left(r,t\right) = -\frac{r}{t} ,
\end{equation}
and resembles some of the shock-free flow scenarios disseminated by, for example, Coggeshall~\cite{coggeshall1991analytic}. With Eq.~(\ref{jnegxi}), Eqs.~(\ref{ode_w_case4}) and (\ref{ode_j_case4}) become, respectively,
\begin{eqnarray}
w' &=& -\frac{\left(n+1\right)\xi w}{2\xi^2}, \label{smooth_w_ode} \\
-1 &=& \frac{K_S w'}{2w^2 \xi}, \label{smooth_j_ode}
\end{eqnarray}
indicating that the problem geometry (represented by $n$) and the EOS (represented by $K_S$) must be mutually constrained to enable a non-trivial solution for $w\left(\xi\right)$. The exact form of this constraint will depend on the exact form the EOS assumes.  

To further illustrate this phenomenon, an example isentropic EOS is the modified Tait EOS given by Eq.~(\ref{tait_eos}). This results in, with Eq. (\ref{bulkmod}),
\begin{eqnarray}
K_S &=& B \gamma \left(\frac{\rho}{\rho_{\rm{ref}}}\right)^\gamma \\
    &=& \frac{B \gamma}{\rho_{\rm{ref}}^{\gamma}} w^{\gamma} ,
\end{eqnarray}
which may then be substituted into Eq.~(\ref{smooth_j_ode}) to yield
\begin{eqnarray}\label{smooth_j_tait_ode}
-1 &=& \frac{\omega w^{\gamma-2} w'}{2\xi} ,
\end{eqnarray}
where 
\begin{equation}\label{omega_def}
\omega = \frac{B \gamma}{\rho_{\rm{ref}}^{\gamma}} .
\end{equation}
Otherwise, Eq.~(\ref{smooth_w_ode}) remains unchanged. The solutions of Eqs.~(\ref{smooth_w_ode}) and (\ref{smooth_j_tait_ode}) are, respectively, 
\begin{eqnarray}
w\left(\xi\right) &=& w_0 \xi^{-\frac{n+1}{2}} , \\
w\left(\xi\right) &=& \left[\left(1-\gamma\right)\left(\frac{\xi^2}{\omega}-w_1\right)\right]^{\frac{1}{\gamma-1}} ,
\end{eqnarray}
where $w_0$ and $w_1$ are arbitrary constants of integration. Simultaneous satisfaction of both forms of $w$ then requires 
\begin{eqnarray}
w_1 &=& 0, \\
w_0 &=& \left(\frac{1-\gamma}{\omega}\right)^{\frac{1}{\gamma-1}},\label{C1_def}\\
\gamma &=& \frac{n-3}{n+1},\label{gamma_n_shock_free}
\end{eqnarray}
thus yielding a constraint involving both $n$ and the material-dependent parameter $\gamma$. Given that the problem geometry factor $n$ may only assume the values $n = 0, 1$, or $2$, only specific modified Tait EOS forms give rise to shock-free, homentropic flows featuring density fields of the form given by Eq.~(\ref{speed_smooth_condition}). Equation (\ref{gamma_n_shock_free}) indicates that $\gamma <0$ in all cases. With Eq.~(\ref{omega_def}), this constraint then indicates $B < 0$ for $\omega$ to be positive definite, and thus $w_0$ given by Eq.~(\ref{C1_def}) to be real-valued. 

When all of above conditions are satisfied, for Case~IV $\xi = \frac{r}{t}$ and $w\left(\xi\right) = \rho\left(r,t\right)$; moreover $P\left(r,t\right)$ and $I\left(r,t\right)$ are given by Eqs.~(\ref{tait_eos}) and (\ref{Isolve}), respectively. The solution in this case is thus comprised of $u\left(r,t\right)$ given by Eq.~(\ref{speed_smooth_condition}), and
\begin{eqnarray}
\rho(r,t) &=& w_0\left(\frac{r}{t}\right)^{\frac{2}{\gamma-1}}, \label{rho_final_smooth}\\
P\left(r,t\right)&=& B\left[\left(\frac{w_0}{\rho_{\rm{ref}}}\right) ^{\gamma}\left(\frac{r}{t}\right)^{\frac{2\gamma}{\gamma-1}} -1\right],\label{P_final_smooth}\\
I(r,t) &=& I_0 + \frac{B}{\rho_{\rm{ref}}} 
\left[\frac{w_0}{\rho_{\rm{ref}}} \left(\frac{r}{t}\right)^{\frac{2}{\gamma-1}}\right]^{-1}\nonumber\\
&&\times\left\{1+\frac{1}{\gamma -1}\left[\left(\frac{w_0}{\rho_{\rm{ref}}}\right)^\gamma \left(\frac{r}{t}\right)^{\frac{2\gamma}{\gamma-1}}\right]\right\}.\label{sie_final_smooth}\;\;\;\;\;
\end{eqnarray}

For the notional parameterization $\rho_{\rm{ref}} = 1.0~\frac{\rm{g}}{\rm{cm^3}}$ and $B = -10^{-3}~\rm{Mbar}$, Eqs.~(\ref{speed_smooth_condition}) and (\ref{rho_final_smooth})-(\ref{sie_final_smooth}) are plotted in Fig.~\ref{smooth_plots} as functions of $\xi$ for all three choices of $n$. Plotting these solutions as functions of $\xi$ clearly reveals the self-similar nature of the flow field: the various shapes depicted in Fig.~\ref{smooth_plots} will hold for any choice of $r$ and $t$ (except $t = 0$, when the flow field is unbounded), aside from a change of scale. Furthermore, as previously discussed for this solution, the geometry factor $n$ sets a unique choice of $\gamma$ in the modified Tait EOS, thus controlling the shape of the density, pressure, and SIE profiles in each case. 

Figure~\ref{smooth_plots} also depicts the pressure field reaching a zero value at a different $\xi$-position for each value of $n$ (or $\gamma$); beyond these points the solution ceases to have physical meaning. The flow field given by Eqs.~(\ref{speed_smooth_condition}) and (\ref{rho_final_smooth})-(\ref{sie_final_smooth}) may thus be interpreted to terminate when $P=0$. With Eq.~(\ref{P_final_smooth}), the $\xi$-position $\xi_0$ where this phenomenon occurs is given by
\begin{equation}
\xi_0 = \left(\frac{\rho_{\rm{ref}}}{w_0}\right)^{\frac{\gamma-1}{2\gamma}} .
\end{equation}
With the Case~IV definition of $\xi = \frac{r}{t}$, $\xi_0$ may be alternatively realized as the space-time trajectory of the zero-pressure surface $r_0\left(t\right)$:
\begin{equation}
r_0\left(t\right) = \xi_0 t .
\end{equation}
Since by construction $\xi_0 > 0$, Eqs.~(\ref{speed_smooth_condition}) and (\ref{rho_final_smooth})-(\ref{sie_final_smooth}) terminated at $r = r_0\left(t\right)$ may thus be interpreted as an expanding bubble solution. Inside of the bubble, the flow field is constrianed to obey Eqs.~(\ref{speed_smooth_condition}) and (\ref{rho_final_smooth})-(\ref{sie_final_smooth}).

%%%%%%%%%%%%%%%%%%%%%%%%%%%%%%%%%%%%%%%%%%
%					 %
%       The Classical Noh Problem 	 %
%					 %
%%%%%%%%%%%%%%%%%%%%%%%%%%%%%%%%%%%%%%%%%%
\subsection{The Classical Noh Problem}\label{noh}
\begin{figure*}[ht]
\includegraphics[scale=0.4]{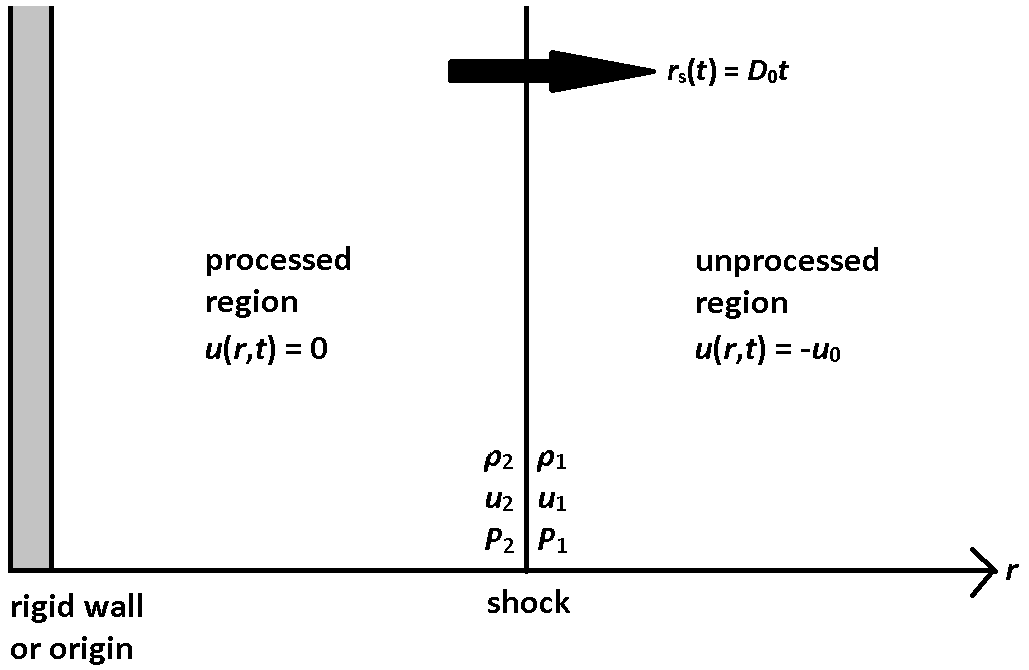}
\caption{Notional depiction of the Noh problem.}
\label{noh_fig}
\end{figure*}
First introduced by W. Noh in 1987~\cite{noh1987errors}, the Noh problem has become ``the workhorse of compressible hydrocode verification for over three decades"~\cite{velikovich2018solution}. Its distinguishing features, potential uses, advantages, disadvantages, physical implications, connections to other physical scenarios, possible generalizations, and a variety of related topics have been extensively documented; see Ramsey et al.~\cite{Ramsey_Boyd_Burnett}, Velikovich and Giuliani~\cite{velikovich2018solution}, and references therein for additional details. 

Of principal interest to this work is the Noh problem formulation as depicted in Fig.~\ref{noh_fig}. The distinguishing features of this scenario are:

\begin{itemize}
\item A constant inflow velocity for all times,
\item An initially constant inflow density. 
\end{itemize}
A less restrictive variant of the Noh problem allowing the inflow velocity to vary in space and time has recently been investigated by Velikovich and Giuliani~\cite{velikovich2018solution}, so we will refer to the traditional scenario as the ``classical" Noh problem.

In this scenario, the impingement of the inflow against the rigid wall at $r = 0$ (or origin in 1D curvlinear geometries) gives rise to a shock wave that propagates outward into the still-incoming fluid. This shock wave brings the fluid behind it to rest. If the impingement at $r = 0$ occurs at $t = 0$, for $t > 0$ the global velocity field $u\left(r,t\right)$ may be written as
\begin{eqnarray}
u(r,t)&=&u_1=-u_0,\;\;\;\;\;\;\;\;\;r_{\rm{s}}<r<\infty, \label{unshocked_speed}\\
u(r,t)&=&u_2=0, \;\;\;\;\;\;\;\;\;\;\;\;\;\; 0<r<r_{\rm{s}}, \label{shocked_speed}
\end{eqnarray}
where the subscripts 1 and 2 again denote the unshocked and shocked regions separated by the time-dependent shock position $r_{\rm{s}}\left(t\right)$, and $u_0$ is a positive constant. The only other constraint featured in the problem formulation is that given by the second distinguishing feature noted above, i.e.,
\begin{equation}
\rho\left(r,t=0\right) = \rho_0, \label{constant_density}
\end{equation}
where $\rho_0$ is a positive constant. 

Conditions for the existence of semi-analytic or even closed-form solutions to the classical Noh problem in any of the 1D geometries (i.e., $n =$ 0, 1, or 2) have been discussed at length by Axford~\cite{axford}, Ramsey et al.~\cite{Ramsey_Boyd_Burnett}, Burnett et al.~\cite{burnett2018verification}, and Velikovich and Giuliani~\cite{velikovich2018solution}. The existence of the 1D planar ($n = 0$) solutions is a direct consequence of arguments advanced by Courant and Friedrichs~\cite{courant1999supersonic} and Menikoff and Plohr~\cite{menikoff1989riemann} for generalized Riemann problems (including shock-piston problems as a special case), though these arguments can also be cast in terms of the universal symmetries inherent Eqs.~(\ref{cons_mass})-(\ref{cons_energy_e_n_p}). Similar conditions in 1D cylindrical or spherical geometries ($n = 1$ or 2, respectively) are more complicated, and require tighter constraints on the EOS closure models coupled to Eqs.~(\ref{cons_mass})-(\ref{cons_energy_e_n_p}).        

Existing mathematical arguments regarding the matter are further illuminated by an intuitive physical interpretation: in the curvilinear geometries, the constant velocity inflow within the unshocked region carries fluid parcels of constant mass into progressively smaller volumes. As time evolves, the fluid density in the unshocked region must therefore increase with decreasing $r$. This variable density field ostensibly gives rise to a variable pressure field, which through momentum conservation would invalidate the constant velocity assumption. The only means through which a self-consistent curvilinear solution can be restored are:

\begin{itemize}
\item Relax the constant velocity assumption inherent to the classical Noh problem, as done by Velikovich and Giuliani~\cite{velikovich2018solution},
\item Employ an EOS allowing for a simultaneous constant pressure and non-constant density state. Several examples of EOS closure laws featuring this property are provided by Axford~\cite{axford}, Ramsey et al~\cite{Ramsey_Boyd_Burnett}, and Burnett et al~\cite{burnett2018verification}.
\end{itemize}

To temporarily obviate these complications in proceeding with the construction of an example solution, we restrict our attention to 1D planar geometry. As discussed in Sec.~\ref{scaling_analysis}, for a classical Noh problem to exist as a manifestation of the scaling group represented by Eq.~(\ref{final_final_group_generator}), any conditions particular to it (i.e., Eqs.~(\ref{shocked_speed}) and (\ref{constant_density})) must be invariant under the action of the group generator $\chi$. To begin, we can rewrite Eqs.~(\ref{unshocked_speed}) and (\ref{constant_density}) as
\begin{eqnarray}
\Upsilon_1 &\equiv& u+u_0 = 0, \label{upsilon1} \\
\Upsilon_2 &\equiv& \rho-\rho_0=0. \label{upsilon2}
\end{eqnarray}
Applying Eq.~(\ref{final_final_group_generator}) to Eq.~(\ref{upsilon1}) (when Eq.~(\ref{upsilon1}) is itself satisfied) yields
\begin{equation}
(a_2-a_1)u_0=0.
\end{equation}
This condition can only be satisfied nontrivially when $a_2=a_1$. It follows that applying Eq.~(\ref{final_final_group_generator}) to Eq.~(\ref{upsilon2}) (when Eq.~(\ref{upsilon2}) is itself satisfied) yields
\begin{equation}
a_3\rho_0=0.
\end{equation}
The only nontrivial solution of this condition is $a_3=0$. 

As a result of the invariance of the classical 1D planar Noh problem's distinguishing conditions, Eq.~(\ref{final_final_group_generator}) reduces to
\begin{eqnarray}\label{noh_group_generator}
\chi=a_1 t \frac{\partial}{\partial t} + a_1 r \frac{\partial}{\partial r},
\end{eqnarray}
corresponding to Case~4 as disseminated in Secs.~\ref{case_4} and \ref{case_4_reduction}. The fact that an arbitrary $K_S$ (and thus EOS) is admissible for the classical 1D planar Noh problem is thus immediately evident from symmetry considerations. The similarity variables associated with this scenario are given by Eqs.~(\ref{xi_def})-(\ref{I_sim_sol}) with $\zeta = 0$, $\beta = 0$, and $\alpha = 1$, and the relevant ODEs for use in the construction of a piecewise solution are given by Eqs.~(\ref{ode_w_case4}) and (\ref{ode_j_case4}). 

%
%
%%%%%%%%%%%%%%%%%%%%%%%%%%%%%%%%%%%%%%%%%%
%					 %
%	      Unshocked region	 	 %
%					 %
%%%%%%%%%%%%%%%%%%%%%%%%%%%%%%%%%%%%%%%%%%
\subsubsection{The Unshocked Region}\label{sec:noh_unshocked}

For the classical Noh problem, the unshocked region is constrained to obey Eq.~(\ref{unshocked_speed}). With this condition $j\left(\xi\right) = -u_0$, and Eq.~(\ref{ode_w_case4}) becomes
\begin{equation}\label{w_noh_unshocked_ode}
w'= 0,
\end{equation}
the solution of which is $w\left(\xi\right) = \rho\left(r,t\right) = \rm{const}$. Given Eq.~(\ref{constant_density}), this solution yields $\rho\left(r,t\right) = \rho_0$ throughout the unshocked region. With this result, Eq.~(\ref{ode_j_case4}) is then identically satisfied. 

Finally, the pressure and SIE in the unshocked region are then given by Eqs.~(\ref{eos}) and (\ref{Isolve}) with $\rho = \rho_0$. Like the fluid density and velocity, these properties are constant throughout the unshocked region. 

%%%%%%%%%%%%%%%%%%%%%%%%%%%%%%%%%%%%%%%%%%
%					 %
%	      Shocked region	 	 %
%					 %
%%%%%%%%%%%%%%%%%%%%%%%%%%%%%%%%%%%%%%%%%%
\subsubsection{The Shocked Region}\label{sec:noh_shocked}

For the classical Noh problem, the shocked region is constrained to obey Eq.~(\ref{shocked_speed}). With this condition $j\left(\xi\right) = 0$, and Eq.~(\ref{ode_w_case4}) again becomes Eq.~(\ref{w_noh_unshocked_ode}). The solution of this equation is again $w\left(\xi\right) = \rho\left(r,t\right) = \rm{const}$. However, unlike the unshocked region, there is no initial condition on the fluid density of the shocked region. As such, the density throughout this region will be denoted $w\left(\xi\right) = \rho\left(r,t\right) = \rho_2$, where $\rho_2$ is a constant to be determined. With this result, Eq.~(\ref{ode_j_case4}) is once again identically satisfied. 

As before, the pressure and SIE in the shocked region are then given by Eqs.~(\ref{eos}) and (\ref{Isolve}) with $\rho = \rho_2$. Like the fluid density and velocity, these properties are constant throughout the shocked region.

%
%
%%%%%%%%%%%%%%%%%%%%%%%%%%%%%%%%%%%%%%%%%%
%					 %
%	      Jump conditions	 	 %
%					 %
%%%%%%%%%%%%%%%%%%%%%%%%%%%%%%%%%%%%%%%%%%
\subsubsection{Rankine-Hugoniot Jump Conditions}\label{sec:noh_jump}

Symmetry analysis of the Rankine-Hugoniot jump conditions provided in Sec.~\ref{jump} yields Eq.~(\ref{D_solution}) as the scale invariant shock speed. With $a_2 = a_1$ for the classical 1D planar Noh problem, Eq.~(\ref{D_solution}) becomes
\begin{equation}
D\left(t\right) = D_0,
\end{equation}
i.e., the classical 1D planar Noh problem features a constant shock velocity (which is to be determined). This outcome is also intuitive on physical grounds, as the constant inflow velocity in the unshocked region must give rise to a constant speed stagnation shock.

With this result, Eq.~(\ref{shocked_speed}), knowledge of the constant state throughout the entire unshocked region, and the isentropic EOS given by Eq.~(\ref{eos}), Eqs.~(\ref{rankine_hugoniot_mass_jump}) and (\ref{rankine_hugoniot_mom_jump}) themselves become two algebraic equations in the two unknowns given by the post-shock density $\rho_2$ and constant shock speed $D_0$:
\begin{eqnarray}
D_0 &=& \frac{\rho_0 u_0}{\rho_2 - \rho_0}, \label{final_noh_mass_hug}\\
f(\rho_2) &=& f(\rho_0)+\rho_0\left(u_0+D_0\right)u_0. \label{final_noh_mom_hug}
\end{eqnarray}
Inserting Eq.~(\ref{final_noh_mass_hug}) into Eq.~(\ref{final_noh_mom_hug}) then yields an algebraic equation exclusively in terms of $\rho_2$:
\begin{equation} \label{transcendental_f}
f(\rho_2) = f(\rho_0)+\rho_0 u_{0}^{2} \left( \frac{1}{1 - \frac{\rho_0}{\rho_2}} \right) ,
\end{equation}
which, depending on the form of the isentropic EOS $f\left(\rho\right)$, is potentially a transcendental algebraic equation for $\rho_2$. With a solution to this equation, the shock velocity may be computed via Eq.~(\ref{final_noh_mass_hug}), the post-shock pressure via Eq.~(\ref{eos}), and the post-shock SIE via Eq.~(\ref{Isolve}). While this solution provides only the immediate post-shock state, given the developments of Sec.~\ref{sec:noh_shocked} it also represents the constant state of the entire shocked region for the classical Noh problem. 
%%%%%%%%%%%%%%%%%%%%%%%%%%%%%%%%%%%%%%%%%%
%					 %
%	      Tait Example	 	 %
%					 %
%%%%%%%%%%%%%%%%%%%%%%%%%%%%%%%%%%%%%%%%%%
\subsubsection{Modified Tait EOS and the Noh Problem}
%\onecolumngrid

%\setlength{\extrarowheight}{1.0ex}
\begin{figure*}[]
	\begin{subfigure}[t]{0.4\textwidth}
	%\centering
        \includegraphics[scale=0.45]{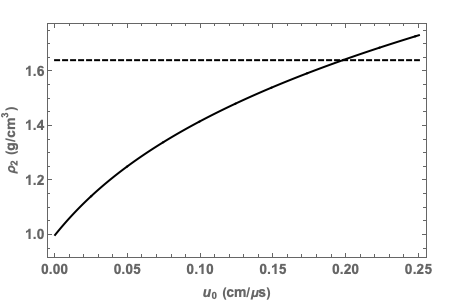}
        \label{density_water_plot}
    \end{subfigure}
        \begin{subfigure}[t]{0.4\textwidth}
 %       \centering
        \includegraphics[scale=0.45]{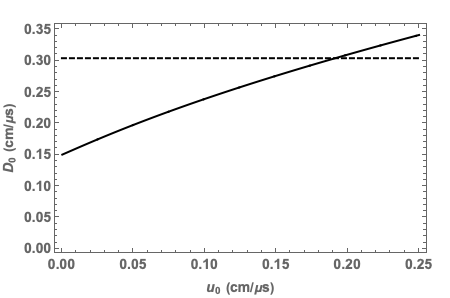}
        \label{speed_water_plot}
    \end{subfigure}   
    \begin{subfigure}[t]{0.4\textwidth}
 %   \centering
        \includegraphics[scale=0.45]{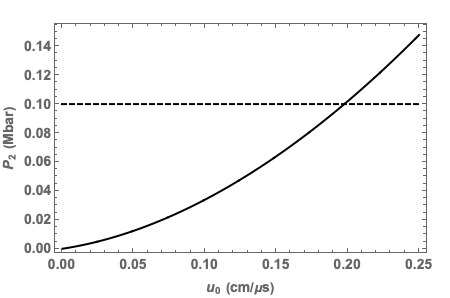}
        \label{pressure_water_plot}
    \end{subfigure}
    \begin{subfigure}[t]{0.4\textwidth}
 %   \centering
        \includegraphics[scale=0.45]{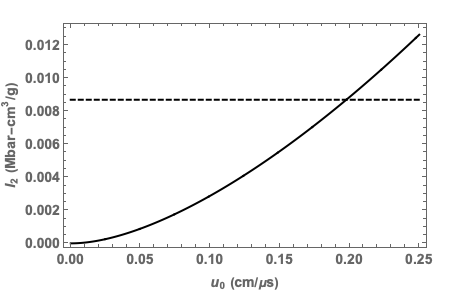}
        \label{sie_water_plot}
    \end{subfigure}
\caption{Shocked state of the classical 1D planar Noh problem for the modified Tait EOS. For a given value of $u_0$, the indicated value of each state variable holds throughout the entire shocked region, and the constant shock speed is as indicated. The dashed horizontal line near the top of each panel corresponds to the maximum pressure at which the modified Tait EOS is assumed to be valid (see Zel'dovich and Raizer~\cite{zeldovich_raizer}). The figures show the shocked density ($\rho_2$), speed ($D_0$), pressure ($P_2$), and SIE ($I_2$) as a function of inflow velocity ($u_0$).}
\end{figure*}
%\twocolumngrid

As discussed in Sec.~\ref{intro}, a canonical example of an isentropic EOS is the modified Tait EOS given by Eq.~(\ref{tait_eos}). The modified Tait form for $f(\rho)$ given by Eq.~(\ref{tait_eos}) may be substituted into all of the otherwise arbitrary developments of Sec.~\ref{noh} to construct an example solution for the classical Noh problem. As discussed in Sec.~\ref{noh}, the construction of a 1D planar instantitation of this problem for an isentropic EOS is piecewise constant, and is essentially encapsulated in a $\rho_2 > \rho_0$ solution of Eq.~(\ref{transcendental_f}). For the modified Tait EOS, this relation becomes
\begin{equation} \label{transcendental_tait}
B\left[\left(\frac{\rho_2}{\rho_{\rm{ref}}}\right)^\gamma - \left(\frac{\rho_0}{\rho_{\rm{ref}}}\right)^\gamma \right] = \rho_0 u_{0}^{2} \left( \frac{1}{1 - \frac{\rho_0}{\rho_2}} \right) ,
\end{equation}
which, given numerical values for the inflow velocity $u_0$, unshocked density $\rho_0$, and the material-dependent constants $\rho_{\rm{ref}}$, $B$, and $\gamma$, can be solved using a root extraction technique. 

As an example, Zel'dovich and Raizer~\cite{zeldovich_raizer} give the following modified Tait parameterization for water:
\begin{eqnarray}
\rho_{\rm{ref}} &=& 1.0 \rm{\frac{g}{cm^3}},\nonumber\\
B &=& 3.214 \times 10^{-3}\;\rm{Mbar},\nonumber\\
\gamma &=& 7.0.\nonumber
\end{eqnarray}
Assuming $\rho = \rho_0 = \rho_{\rm{ref}}$ throughout the unshocked region, the resulting numerical solution of Eq.~(\ref{transcendental_tait}) is given in Fig.~(\ref{density_water_plot}) for a range of $u_0$ values. The associated shock speed and shocked fluid pressure, calculated using Eqs.~(\ref{final_noh_mass_hug}) and (\ref{tait_eos}), are depicted in Figs.~(\ref{speed_water_plot}) and (\ref{pressure_water_plot}), respectively. The shocked fluid SIE is calculated using Eqs.~(\ref{tait_eos}) and (\ref{Isolve}), which yield
\begin{equation}\label{sie_tait}
I = \frac{f\left(\rho_2\right) + \gamma B}{\left(\gamma -1 \right)\rho_2} + I_0,
\end{equation}
and is depicted in Fig.~\ref{sie_water_plot} for $I_0 = -\frac{\gamma B}{\left(\gamma -1\right)\rho_0}$, so that $I\left(\rho_0\right) = 0$.

Analysis of Eq.~(\ref{transcendental_tait}) shows that a $\rho_2 > \rho_0$ solution is obtained for any $u_0 > 0$. When $u_0$ satisfies this condition, $\rho_2$ increases from $\rho_0$ monotonically and without limit with increasing $u_0$, as depicted in Fig.~\ref{density_water_plot}. The shock velocity likewise increases monotonically and without limit with increasing $u_0$, but as $u_0 \to 0$ it limits to the sound speed $c$ in the unshocked region, which may be calculated using Eqs.~(\ref{tait_eos}), (\ref{bulkmod}), (\ref{sound}), and $\rho = \rho_{\rm{ref}}$ as:
\begin{eqnarray}
c &=& \sqrt{\frac{K_S}{\rho}}, \nonumber \\
  &=& \sqrt{\frac{\gamma B}{\rho_{\rm{ref}}}}. \label{sound_tait}
\end{eqnarray}
The shocked pressure and shocked SIE likewise increase from zero monotonically and without limit with increasing $u_0$.  

The reason for these behaviors is associated with the modified Tait EOS entropy. Using Eqs.~(\ref{tait_eos}), (\ref{bulk_def}), and (\ref{bulkmod}), the entropy in the shocked region may be calculated as
\begin{eqnarray}
S &=& F\left[ \left( P+B\right) \rho^{-\gamma} \right] +\rm{const.}, \nonumber \\
  &=& F\left( B \rho_{\rm{ref}}^{\gamma} \right) +\rm{const.}, \label{tait_entropy1}
\end{eqnarray}
which is explicitly constant for any function $F$ of the indicated argument. In particular, the entropy in the shocked region assumes this value for any inflow velocity $u_0$, and thus is independent of shock strength. Equation~(\ref{tait_entropy1}) thus indicates that given a constant entropy in the shocked region, increasing $u_0$ must always be counterbalanced by increasing $\rho_2$. 

This behavior is distinct from that encountered in the classical ideal gas Noh problem, which features a maximum shock compression limit independent of $u_0$. Following from the definition of the ideal gas entropy given by Zel'dovich and Raizer~\cite{zeldovich_raizer}, and the classical 1D planar ideal gas solution disseminated by Ramsey et al.~\cite{Ramsey_Boyd_Burnett}, the entropy in the shocked region for this case is given by
\begin{eqnarray}
S &=& c_v \ln\left(P\rho^{-\gamma}\right) + \rm{const.}, \nonumber \\
  &=& c_v \ln\left[\frac{1}{2}u_{0}^{2} \rho_{0}^{1-\gamma} \left(\gamma+1\right)\left(\frac{\gamma+1}{\gamma-1}\right)^{-\gamma}\right] + \rm{const.},\;\;\;\;\;\;\;\;\;\;\label{ideal_entropy}
\end{eqnarray}
where $c_v$ and $\gamma$ are interpreted as the constant specific heat capacity at constant volume and constant adiabatic index, respectively. This result clearly depends on $u_0$ and is thus not explicitly constant with respect to that parameter. The difference in behavior of $\rho_2\left(u_0\right)$ between these two cases is thus clearly revealed to be one consequence of invoking an isentropic EOS assumption in the construction of a classical 1D planar Noh problem.
%%%%%%%%%%%%%%%%%%%%%%%%%%%%%%%%%%%%%%%%%%
%					 %
%	       Conclusion	 	 %
%					 %
%%%%%%%%%%%%%%%%%%%%%%%%%%%%%%%%%%%%%%%%%%
\section{Discussion and Conclusion}\label{conclusion}

Under the assumption of an isentropic EOS of the form given by Eq.~(\ref{eos}), the 1D inviscid Euler equations may be collapsed to a set of two coupled, nonlinear PDEs. In this formulation total energy conservation is automatically ensured, as is isentropic flow. A special case of this phenomenon involving explicitly constant entropy is referred to as homentropic flow. Piecewise isentropic or homentropic flows featuring shock waves may also exist in these scenarios, with the shocked and unshocked states connected by an appropriate form of the Rankine-Hugoniot jump conditions.

Moreover, the equations governing any of the aforementioned flows may be subjected to symmetry analysis, in the interest of reducing the PDEs and any ancillary conditions to simpler structures (e.g., ODEs) more easily amenable to either exact or semi-analytical solution. For cases where the included isentropic EOS (encoded in an adibatic bulk modulus $K_S$) is left as an arbitrary function of the fluid density $\rho$, symmetry analysis yields conditional forms the EOS may assume so as to ensure the presence of various symmetries (e.g., scaling transformations). When the isentropic EOS assumes one of these forms, the PDEs governing the associated fluid motion may be reduced to ODEs.

The ODEs obtained via symmetry analysis are likely easier to solve than their PDE counterparts, and thus may be used to construct a variety of exact or semi-analytical solutions with desired properties. Under scaling transformations, one such example is the classical Noh problem featuring a constant velocity inflow directed against a rigid wall (1D planar geometry) or curvilinear origin (1D cylindrical or spherical geometries), giving rise to an outward propagating, constant velocity stagnation shock. For the case of an arbitrary isentropic EOS, the solution of this problem essentially reduces to a transcendental solve in the shocked density $\rho_2$, which may then be used to reconstruct the shock trajectory and entire shocked flow field.  

\subsection{Recommendations for Future Work}\label{future_work}

The modified Tait EOS given by Eq.~(\ref{tait_eos}) is similar in form to the ideal gas EOS, namely
\begin{equation}\label{ideal_gas}
P = \left(\gamma -1\right)\rho I ,
\end{equation}
in that the adibatic bulk moduli calculated from them are only slightly different in form. A potentially fruitful avenue for future work would be to further assess the consequences these differences manifest in solutions of various Noh-like problems, and to investigate the conditions under which one of the result sets can be obtained from the other. Similar efforts could be performed with respect to the recent work of both Deschner et al.~\cite{deschner2018self} and Velikovich and Giuliani~\cite{velikovich2018solution}.  

Additional natural extensions of this work include (but are not limited to):

\begin{itemize}
\item Extension of the current results for the modified Tait EOS to other, similar EOS examples with validity in a variety of regimes and contexts. One such example is given by the Birch-Murnaghan EOS, as reported by Birch~\cite{Birch}.
\item Determination of the conditions for scale-invariance of the homentropic Euler equations in other coordinate systems (e.g., 2D or 3D). Such an analysis should follow easily from that performed in Sec.~\ref{scaling_analysis}, and can also be connected to various outcomes reported by Ovsiannkov~\cite{ovsyannikov1} or Holm~\cite{holm1976symmetry}.
\item Determination of conditions for the presence of any symmetries (i.e., not limited to scaling transformations) in the homentropic Euler equations, in any coordinate system. For example, in addition to the kinematic scaling ($r-t$) transformation indicated by Eq.~(\ref{noh_group_generator}), the 1D planar instantiation of Eqs.~(\ref{simp_eos_cons_mass}) and (\ref{simp_eos_cons_mom}) is for any EOS invariant under time translation, space translation, and Galilean boost transformations (see, for example, Axford~\cite{axford}). The association with and interplay of these various transformations with canonical solutions remains to be rigorously assessed. 
\item Construction of Sedov, Guderley, or other analogous test problems featuring shock waves in the piecewise homentropic setting, following from the presence of any symmetries.
\item Construction of additional shock-free solutions (e.g., in the style of Coggeshall~\cite{coggeshall1991analytic} or McHardy et al.~\cite{mchardy2019group}) of the homentropic Euler equations in any coordinate system.
\item Establishment of the connections between any of the symmetries or associated new solutions described above, and their ideal gas counterparts (if they exist).
\end{itemize}

Given the wide scope of potential work available in the context of the isentropic Euler equations, this work may serve as a foundation for any future developments and applications along these lines.
\section{Acknowledgements}
This work was supported by the U.S. Department of Energy (DOE) through the Los Alamos National Laboratory. Los Alamos National Laboratory is operated by Triad National Security, LLC, for the National Nuclear Security Administration of the DOE (contract number 89233218CNA000001). JFG was partially funded through supported by the National Science Foundation (NSF) under Grant No. PHY-1803912. The authors would like to thank E. J. Albright, B. A. Temple, J. D. McHardy, P. J. Jaegers, E. M. Schmidt, J. H. Schmidt, and J. A. Tellez for their valuable insights on these topics. We would also like to thank the two reviewers of this article for their insight.
%%%%%%%%%%%%%%%%%%%%%%%%%%%%%%%%%%%%%%%%%%
%					 %
%	       Appendix		 	 %
%					 %
%%%%%%%%%%%%%%%%%%%%%%%%%%%%%%%%%%%%%%%%%%
\appendix
\section{Similarity Variables Substitution into the Euler Equations}\label{sim_sub_euler}
The following relations will be useful in the calculation of Eqs. (\ref{new_rho_w}) and (\ref{new_u_j}) which are: 
\begin{eqnarray}
\beta &=&1-\frac{1}{\alpha},\label{alpha_beta_relationship}
\end{eqnarray}
and
\begin{eqnarray}
t^{-1} &=& \xi^{\frac{1}{\alpha}}r^{-\frac{1}{\alpha}}=\xi^{\frac{1}{\alpha}}r^{\beta-1},\label{inverse_t_relationship}
\end{eqnarray}
and
\begin{eqnarray}
\frac{\partial \xi}{\partial t}&=& -\alpha r t^{-\alpha-1},\\
\frac{\partial \xi}{\partial r}&=& t^{\alpha}.
\end{eqnarray}
Taking each respective time and space derivative of Eqs. (\ref{rho_sim_sol}) and (\ref{u_sim_sol}), we have:
\begin{align}
\frac{\partial\rho}{\partial t} &= r^\zeta \frac{\partial\xi}{\partial t}\frac{dw}{d\xi},\nonumber\\
&=-\alpha r^{\zeta+1}t^{-\alpha-1}\frac{dw}{d\xi},\\
\frac{\partial\rho}{\partial r} &=\zeta r^{\zeta-1}w + r^\zeta \frac{\partial \xi}{\partial r}\frac{dw}{d\xi},\nonumber\\
&=\zeta r^{\zeta-1}w + r^\zeta t^{-\alpha}\frac{dw}{d\xi},\\
\frac{\partial u}{\partial t} &=r^\beta\frac{\partial \xi}{\partial t}\frac{dj}{d\xi},\nonumber\\
&=-\alpha r^{\beta+1}t^{-\alpha-1}\frac{dj}{d\xi},\\
\frac{\partial u}{\partial r} &=\beta r^{\beta-1}j + r^\beta\frac{\partial \xi}{\partial r}\frac{dj}{d\xi},\nonumber\\
&=\beta r^{\beta-1}j + r^\beta t^{-\alpha}\frac{dj}{d\xi}. 
\end{align}
We first perform the substitution for Eq. (\ref{simp_eos_cons_mass}), that is
\begin{widetext}
\begin{eqnarray}
0&=&r\frac{\partial \rho}{\partial t}+ru\frac{\partial \rho}{\partial r}+r\rho\frac{\partial u}{\partial r} + n\rho u,\\
0&=&r\left[ -\alpha r^{\zeta+1}t^{-\alpha-1}\frac{dw}{d\xi}\right]+rr^{\beta}j\left[\zeta r^{\zeta-1}w + r^\zeta t^{-\alpha}\frac{dw}{d\xi} \right] + r r^\zeta w \left[\beta r^{\beta-1}j + r^\beta t^{-\alpha}\frac{dj}{d\xi}\right]\nonumber\\&&+nr^\zeta r^\beta w j\\
0&=&-\alpha r^{\zeta+1}t^{-1}\xi \frac{dw}{d\xi}+\zeta r^{\zeta+\beta}jw + r^{\zeta+\beta}\xi j\frac{dw}{d\xi}+\beta r^{\zeta+\beta}jw + r^{\zeta+\beta}\xi w \frac{dj}{d\xi}+n r^{\zeta+\beta}wj,
\end{eqnarray}
\end{widetext}
using Eq. (\ref{inverse_t_relationship}) we have
\pagebreak
\begin{widetext}
\begin{eqnarray}
0&=&-\alpha r^{\zeta+\beta}\xi^{1+\frac{1}{\alpha}}\frac{dw}{d\xi} + \left[\zeta+\beta+n\right]r^{\zeta+\beta}jw + r^{\zeta+\beta}\xi\left[w\frac{dj}{d\xi}+j\frac{dw}{d\xi}\right],\\
0&=&\left[-\alpha\xi^{1+\frac{1}{\alpha}}+\xi j\right]\frac{dw}{d\xi} + \left[\zeta+\beta+n\right]jw + \xi w\frac{dj}{d\xi},\\
\frac{dw}{d\xi}&=&\frac{\xi w(\xi)\frac{dj}{d\xi}+\left(\zeta+\beta+n\right)j(\xi)w(\xi)}{\xi\left(\alpha\xi^{\frac{1}{\alpha}}-j\right)},
\end{eqnarray}
\end{widetext}
which is Eq. (\ref{new_rho_w}).

Finally we perform the substitution for Eq. (\ref{simp_eos_cons_mom}). It follows,
\begin{widetext}
\begin{eqnarray}
0&=&\rho^2\frac{\partial u}{\partial t} + \rho^2 u\frac{\partial u}{\partial r} +K_S(\rho)\frac{\partial\rho}{\partial r},\\
0&=&\rho^{2\zeta}w^2\left[-\alpha r^{\beta+1}t^{-\alpha-1}\frac{dj}{d\xi}\right]+r^{2\zeta}r^{\beta}jw^2\left[\beta r^{\beta-1}j+r^\beta t^{-\alpha}\frac{dj}{d\xi}\right]+K_S[r^\zeta w]\bigg[\zeta r^{\zeta-1}w+r^\zeta t^{-\alpha}\frac{dw}{d\xi}\bigg],\\
0&=&-\alpha r^{2\zeta+\beta}t^{-1}\xi w^2\frac{dj}{d\xi}+\beta r^{2\zeta+2\beta-1}w^2j^2+r^{2\zeta+2\beta}t^{-\alpha} w^2 j\frac{dj}{d\xi}+\zeta K_S[r^\zeta w]r^{\zeta-1}w+K_s(r^\zeta w)r^{\zeta-1}\xi\frac{dw}{d\xi},
\end{eqnarray}
\end{widetext}
using Eq. (\ref{inverse_t_relationship}) we have
\begin{widetext}
\begin{eqnarray}
0&=&-\alpha r^{2\zeta+2\beta-1}\xi^{1+\frac{1}{\alpha}}w^2\frac{dj}{d\xi}+\beta r^{2\zeta+2\beta-1}w^2j^2+r^{2\zeta+2\beta-1}\xi w^2 j \frac{dj}{d\xi}+\zeta K_S[r^\zeta w]r^{\zeta-1}w+K_S[r^\zeta w]r^{\zeta-1}\xi \frac{dw}{d\xi},\\
0&=&\left[-\alpha w^2\xi^{1+\frac{1}{\alpha}} + \xi w^2\right]\frac{dj}{d\xi} + \beta w^2j^2+r^{-\zeta-2\beta}K_S[r^\zeta w]\left(\zeta w +\xi\frac{dw}{d\xi}\right),\\
\frac{dj}{d\xi}&=&\frac{r^{-\zeta-2\beta}K_S\left[r^\zeta w(\xi)\right]\left(\zeta w +\xi\frac{dw}{d\xi}\right)+\beta w^2j^2 }{\xi w^2 \left(\alpha\xi^{\frac{1}{\alpha}}-j\right) },
\end{eqnarray}
\end{widetext}
which is Eq. (\ref{new_u_j}).
\bibliographystyle{apsrev4-1}
\bibliography{Isentropic_Euler}{}
\end{document}